\documentclass[final,12pt,times]{elsarticle}

\usepackage[T1]{fontenc}
\usepackage{amsmath,amssymb,amsthm}
\usepackage{graphicx}      
\usepackage{color}
\usepackage{enumerate}
\usepackage{tikz}
\usetikzlibrary{arrows,automata}
\usepackage{enumitem}

\usepackage{lineno}
\usepackage{hyperref}
\usepackage{setspace}
\newtheorem{thm}{\textbf{Theorem}}[section]

\newtheorem{definition}[thm]{\textbf{Definition}}
\newtheorem{lem}[thm]{\textbf{Lemma}}

\newtheorem{prop}[thm]{\textbf{Proposition}}

\newtheorem{remark}[thm]{\textbf{Remark}}

\numberwithin{equation}{section}

\allowdisplaybreaks

\newcommand{\R}{\mathbb{R}}

\newcommand{\C}{\mathbb{C}}

\newcommand{\E}{\mathbb{E}}

\renewcommand{\H}{\mathcal{H}}

\newcommand{\Imax}{I_{\text{max}}}
\newcommand{\x}{\langle x\rangle}

\newcommand{\EE}{\mathbf{E}}

\let\Re=\undefined\DeclareMathOperator*{\Re}{Re}
\let\Im=\undefined\DeclareMathOperator*{\Im}{Im}

\begin{document}

\begin{frontmatter}
	\journal{Journal of Mathematical Physics}
	
	\title{Almost sure scattering for the energy-critical Hartree equation with radial data below the energy space}
	
	\author[nn]{Yuli Feng}
	\author[nm]{Yuyang Ye}
	\author[nn]{Han Yang\corref{cor1}}
	
	\cortext[cor1]{Corresponding author.
		Tel.: +86 028 66368172.}
	\ead{hanyang95@263.net}
	
	\address[nn]{School of Mathematics, Southwest Jiaotong University, Chengdu 611756, PR China}
	\address[nm]{School of Mathematical Sciences, Fudan University, Shanghai 200433, PR China}
\begin{abstract}
We study the energy-critical Hartree equation on $\mathbb{R}^d$, $d\ge5$, with randomized radial initial data of supercritical regularity. To handle the nonlocal Hartree interaction, we exploit the Riesz smoothing of the Hartree convolution through a bilinear decomposition adapted to the regularity of the interacting terms. This leads to a Hartree bilinear estimate that allows us to combine probabilistic weighted spacetime estimates with weighted and interaction Morawetz estimates. For every $d\ge5$, we establish almost sure local well-posedness for randomized radial data below the energy regularity, together with a conditional scattering criterion. For $d\ge7$, we further prove almost sure global existence and scattering. The uniform energy and Morawetz bounds are established through a double bootstrap argument combining the Hartree bilinear decomposition with an energy increment estimate.
\end{abstract}


\end{frontmatter}
\section{Introduction}

\par\quad

In this paper, we consider the following Cauchy problem for the Hartree equation in dimensions $d \geq 5$:
\begin{equation}
	\label{HE1}
	\begin{cases}
		i\partial_t u + \Delta u = (|x|^{-4} \ast |u|^2)u,\\
		u|_{t = 0} = u_0,
	\end{cases}
	\quad (t,x) \in \mathbb{R} \times \mathbb{R}^d,
\end{equation}
where $\ast$ denotes convolution on $\mathbb{R}^d$.
Hartree equations arise naturally as effective mean-field models for interacting many-body quantum systems and provide a fundamental class of nonlocal nonlinear dispersive equations \citep{Ammari08}.
It is well known that equation \eqref{HE1} enjoys a dilation symmetry. More precisely, if $u(t,x)$ is a solution to \eqref{HE1} with initial condition $u_0$, then
\begin{equation}\label{scaling}
	u_\lambda(t,x)=\lambda^{\frac{d-2}{2}}u(\lambda^2 t,\lambda x), \qquad \lambda>0,
\end{equation}
is also a solution to \eqref{HE1} with the $\lambda$-scaled initial condition. Associated to this dilation symmetry, the $\dot H^s(\mathbb{R}^d)$-norm of initial data behaves as
\begin{equation}\label{H1 norm}
	\|u_\lambda(0,x)\|_{\dot H_x^s(\mathbb{R}^d)}
	=
	\lambda^{s-1}\|u_0\|_{\dot H_x^s(\mathbb{R}^d)}.
\end{equation}
In particular, when $s=1$, the $\dot H_x^1(\mathbb{R}^d)$-norm is invariant under the scaling \eqref{scaling}. For this reason, $\dot H_x^1(\mathbb{R}^d)$ is referred to as the critical Sobolev space for \eqref{HE1}. Moreover, \eqref{HE1} admits the conserved energy
\begin{equation}\label{E}
	E(u(t))
	=\frac12 \|\nabla u(t)\|_{L^2}^2
	+\frac14 \iint_{\mathbb{R}^d \times \mathbb{R}^d}
	\frac{|u(t,x)|^2\,|u(t,y)|^2}{|x-y|^4}\,dx\,dy,
\end{equation}
which is also invariant under the scaling \eqref{scaling}. For these reasons, the Cauchy problem \eqref{HE1} is called energy-critical.

Let us first discuss the known deterministic results. In the literature, the Cauchy problem of the Hartree equation
\begin{equation}\label{NLS normal}
	i\partial_t u + \Delta u = \mu(|x|^{-\gamma} \ast |u|^2)u,
	\quad (t,x) \in \mathbb{R} \times \mathbb{R}^d,
\end{equation}
has been intensively studied, where $\mu = \pm 1$ and $0 < \gamma < d$.
Equation \eqref{NLS normal} is referred to as the defocusing case when 
$\mu = 1$, and as the focusing case when $\mu = -1$. 

Miao, Xu and Zhao
\citep{MiaoXuZhao2008CauchyHartree} established local well-posedness in $H^s(\mathbb{R}^d)$ for
$s\ge\max\{0,\tfrac{\gamma}{2}-1\}$ together with small-data scattering in $H^1$. Scattering in the defocusing
$\dot H^1$-subcritical range was developed by Ginibre and Velo \citep{GinibreVelo2000} following the
method of Morawetz and Strauss \citep{MorawetzCathleenStrauss1972}, and global well-posedness and scattering in
the mass-critical case were obtained in \citep{Miao09}. In the energy-critical case
$\gamma=4$, Miao, Xu and Zhao \citep{Miao07} proved global well-posedness and scattering for
radial $\dot H^1(\mathbb{R}^d)$ data with $d\ge5$, and the radial assumption was
later removed in \citep{Miao11} through a concentration-compactness argument. Further results on \eqref{NLS normal} can be found 
in \citep{Li09, Manna22}.
On the other hand, the deterministic theory for \eqref{HE1} in the 
supercritical regime $s < 1$ is in general ill-posed \citep{MiaoXuZhao2008CauchyHartree}, as is the case 
for power-type nonlinearities \citep{Pausader09}. 
However, probabilistic methods have led to substantial 
progress in constructing both local and global solutions for several nonlinear 
dispersive equations at regularities below the scaling-critical threshold.
Probabilistic methods in dispersive equations go back to Bourgain \citep{Bourgain96}, who
studied random data in connection with the Gibbs measure for the cubic nonlinear Schrödinger equation (NLS) on
$\mathbb{T}^2$. This direction was further developed by Burq and Tzvetkov
\citep{Burq081,Burq082} through their study of random-data Cauchy
theory on compact manifolds. They combined spectral randomization with
invariant measure arguments to establish almost sure local
well-posedness and global existence for supercritical nonlinear wave
equations (NLW). In the two-dimensional radial setting, Deng \citep{Deng12}
established almost sure global well-posedness and invariance of the
Gibbs measure for NLS with a harmonic potential, together with almost
sure global well-posedness and scattering for the corresponding
defocusing equation without the potential.
On $\mathbb R^d$, there is no countable basis of $L^2(\mathbb R^d)$
consisting of eigenfunctions of the Laplacian. B\'enyi, Oh and
Pocovnicu \citep{Benyi15} therefore introduced a randomization adapted
to the Wiener decomposition of frequency space (see
Definition~\ref{def:random}). Using this randomization, they established
almost sure local well-posedness, together with small-data global
well-posedness and scattering, for the cubic NLS below the
scaling-critical regularity \citep{Benyi15,Benyi19}. Almost sure global
existence and well-posedness results for power-type NLW on
$\mathbb R^3$ were obtained by L\"uhrmann and Mendelson
\citep{Lührmann14,Lührmann16}.
More recently, Deng, Nahmod and Yue developed the methods of random
averaging operators and random tensors to study the propagation of
randomness in nonlinear dispersive equations. Their work established
almost sure local well-posedness for semilinear
Schr\"odinger equations throughout the full subcritical range relative
to the probabilistic scaling \citep{DengNahmodYue2022}, as well as almost sure global well-posedness
and invariance of the associated Gibbs measure for Wick-ordered
defocusing NLS on $\mathbb T^2$
\citep{DengNahmodYue2024NLS}. 

In the energy-critical setting, probabilistic global well-posedness has been
studied extensively for NLW and NLS. Pocovnicu \citep{Pocovnicu17} proved almost
sure global existence and uniqueness for energy-critical NLW with rough random
initial data in dimensions $d=4,5$, and Oh and Pocovnicu \citep{Oh16} obtained
almost sure global well-posedness for the defocusing quintic NLW on
$\mathbb{R}^3$ for $s>\frac12$. For the energy-critical NLS, Brereton \citep{Brereton2019} proved almost sure local well-posedness for the quintic NLS on $\mathbb R^3$ below the critical regularity and obtained a conditional global well-posedness criterion. The conditional global arguments in  \citep{Benyi15, Brereton2019} require an a priori bound on the critical Sobolev norm of the nonlinear part of the solution.
Unconditional almost sure global existence and scattering were later
established by Dodson, L\"uhrmann and Mendelson
\citep{Dodson+Luhrmann+Mendelson_2020} for the energy-critical NLW with
randomized radial data, and by Killip, Murphy and Visan
\citep{Killip+Murphy+Visan_2019} for the energy-critical NLS with randomized
radial data. The latter result was further improved in \citep{Dodson19}.
For further probabilistic well-posedness and scattering results for NLS, we
refer to \citep{Thomann09,Poiret14,Murphy19,Zhang2023}.  Motivated by these works, we study almost sure scattering for the energy-critical Hartree equation in dimensions $d\geq5$.
Compared with the local NLS nonlinearity, the Hartree interaction is nonlocal and gives rise to convolution terms whose mixed deterministic-random terms cannot be treated by the standard pointwise product estimates. This nonlocal structure requires bilinear estimates adapted to the Hardy--Littlewood--Sobolev framework and  leads to one of the main difficulties addressed in this work.

Regarding the Hartree equation with random initial data, Chen and Gao
\citep{Chen15} studied the Cauchy problem for Hartree equations in the
energy-subcritical case. Using the Burq--Tzvetkov randomization of the initial
data, they constructed strong solutions below the deterministic critical
threshold. Meanwhile, Deng, Nahmod and Yue
\citep{DengNahmodYue2021Hartree} studied a different probabilistic setting on
$\mathbb T^3$, proving invariance of the Gibbs measure and global existence of
strong solutions. More recently, Tao and Zhao \citep{Tao23} proved almost sure
scattering for the defocusing energy-critical Hartree equation on
$\mathbb{R}^5$, using a ``Narrowed'' Wiener randomization together with
stability theory and a modified interaction Morawetz estimate. In contrast,
we work with the standard unit-scale Wiener randomization introduced in
Definition \ref{def:random}. In this framework, the Hartree potential has Fourier multiplier
$|\xi|^{-(d-4)}$, and the resulting smoothing can be distributed between the
rough random factor and the energy-regular factor. The main difficulty comes from the nonlocal Hartree interaction, whose
convolution structure prevents the pointwise factorization available for local
power-type nonlinearities. We prove almost sure local well-posedness and a
conditional scattering criterion for all $d\ge5$, and almost sure global
existence and scattering for $d\ge7$. The key mechanism in our argument is a bilinear decomposition that distributes the Riesz smoothing associated with the Hartree potential between the two factors. The resulting bilinear estimates allow the probabilistic spacetime bounds to be incorporated into the weighted Morawetz estimate, the interaction Morawetz estimate, and the energy increment estimate. In the present nonlocal setting, they thus play the role of the pointwise product estimates used for local power-type nonlinearities in \citep{Dodson+Luhrmann+Mendelson_2020,Killip+Murphy+Visan_2019}. This mechanism enables us to close a double bootstrap for the energy and Morawetz quantities and thereby obtain the uniform energy bound required for the global argument. Before stating the main results, we introduce the randomization procedure used throughout the paper.

\begin{definition}[Randomization]\label{def:random}
Let $\psi\in C_c^\infty(\mathbb{R}^d)$ be a nonnegative function supported
in $[-1,1]^d$ such that
\[
\sum_{k\in\mathbb{Z}^d}\psi_k(\xi)=1
\qquad\text{for all }\xi\in\mathbb{R}^d,
\]
where $\psi_k(\xi):=\psi(\xi-k)$.
Let $s\in\mathbb{R}$ and $f\in H^s(\mathbb{R}^d)$. For each
$k\in\mathbb{Z}^d$, define
\[
f_k:=\bigl[\widehat f\,\psi_k\bigr]^\vee.
\]
Let $\{X_k\}_{k\in\mathbb{Z}^d}$ be independent standard real-valued
Gaussian random variables defined on a probability space
$(\Omega,\Sigma,\mathbb{P})$. We define the randomization of $f$ by
\[
f^\omega(x)
=
\sum_{k\in\mathbb{Z}^d}X_k(\omega)f_k(x).
\]
\end{definition}
For definiteness, we work with the Gaussian randomization above.
By standard Khintchine-type inequalities, the arguments extend to more
general randomizations satisfying a subgaussian moment condition such as
$
\mathbb{E}\bigl(e^{\gamma X_k}\bigr)
\le e^{c\gamma^2}
$
uniformly in $\gamma\in\mathbb{R}$ and $k\in\mathbb{Z}^d$ for some
constant $c>0$. If $f\in H^s(\mathbb{R}^d)$, then
$f^\omega\in H^s(\mathbb{R}^d)$ almost surely, with
$
\mathbb{E}\|f^\omega\|_{H^s}^2
\sim
\|f\|_{H^s}^2.
$
The randomization is linear and commutes with every Fourier multiplier,
but does not preserve radial symmetry.

Throughout the paper, we write $ u(t) = e^{it\Delta}f^\omega + v(t)$ and $F := e^{it\Delta}f^\omega$. The function $v$ solves the forced equation \begin{equation}\label{forced1} (i\partial_t + \Delta)v = \bigl(|x|^{-4} \ast |v + F|^2\bigr)(v + F) \qquad with \qquad v(0) = 0. \end{equation}
Our main results are the following two theorems.
\begin{thm}[Almost sure local well-posedness and conditional scattering]
	\label{thm:local}
	Let $d\ge 5$ and let $f\in H^s_{\mathrm{rad}}(\mathbb{R}^d)$ with
$
	\frac{d}{2(d-1)}<s<1.
$
	Then for almost every $\omega\in\Omega$, there exists a unique
	maximal-lifespan solution
	\[
	u=e^{it\Delta}f^\omega+v,
	\qquad
	v\in C(I_{\max};H^1(\mathbb{R}^d))
	\cap
	L_{\mathrm{loc}}^4
	\left(
	I_{\max};
	L^{\frac{2d}{d-3}}(\mathbb{R}^d)
	\right),
	\qquad
	0\in I_{\max},
	\]
	to \eqref{HE1} with $u(0)=f^\omega$.
	Moreover, if
$
	\sup_{t\in I_{\max}}\|v(t)\|_{H_x^1}<\infty,
$
	then $I_{\max}=\mathbb{R}$ and the nonlinear part $v$ scatters in
	$H^1(\mathbb{R}^d)$. More precisely, there exist unique
	$v_\pm\in H^1(\mathbb{R}^d)$ such that
	\[
	\lim_{t\to\pm\infty}
	\bigl\|
	u(t)-e^{it\Delta}(f^\omega+v_\pm)
	\bigr\|_{H^1(\mathbb{R}^d)}
	=0.
	\]
\end{thm}

\begin{thm}[Almost sure global existence and scattering]
	\label{thm:global}
	Let $d\ge 7$ and let $f\in H^s_{\mathrm{rad}}(\mathbb{R}^d)$ with
$
	1-\frac{2}{d}<s<1.
$
	Then for almost every $\omega\in\Omega$, there exists a unique global
	solution
	\[
u=e^{it\Delta}f^\omega+v,
\qquad
v\in C(\mathbb R;H^1(\mathbb R^d))
\cap
L_t^4L_x^{\frac{2d}{d-3}}
(\mathbb R\times\mathbb R^d),
	\]
	to \eqref{HE1} with $u(0)=f^\omega$.
	The nonlinear part $v$ scatters in $H^1(\mathbb{R}^d)$. More precisely, there
	exist unique $v_\pm\in H^1(\mathbb{R}^d)$ such that
	\[
	\lim_{t\to\pm\infty}
	\bigl\|u(t)-e^{it\Delta}(f^\omega+v_\pm)\bigr\|_{H^1(\mathbb{R}^d)}
	=0.
	\]
\end{thm}

\begin{remark}
	For $s<1$, the deterministic Cauchy problem \eqref{HE1} is ill-posed
	\citep{MiaoXuZhao2008CauchyHartree}. Thus, in Theorems~\ref{thm:local} and~\ref{thm:global}, we consider
	initial data below the critical regularity $s_c=1$. The two regularity thresholds arise from different probabilistic estimates. The local theory relies only on the weighted gradient estimate in
	Proposition~\ref{P:STE}, which holds for
	$s>\frac{d}{2(d-1)}$. On the other hand, the energy increment requires in addition the $L^\infty_t L^{\frac{2d}{d-2}}_x$ bound of Proposition~\ref{P:STE2}, which holds
	for $s>1-\frac{2}{d}$. 
	In dimensions $d=5,6$, the derivative on the Hartree potential
	$\nabla(|x|^{-4}\ast \cdot)$ introduces a loss of half a derivative
	when $d=6$ and a full derivative when $d=5$. This is the reason for the
	assumption $d\ge7$ in Theorem~\ref{thm:global}. See
	Remark~\ref{rmk:d7} for details.
\end{remark}

\begin{remark}
	Although our argument is motivated by the random-data energy-critical theory
	developed for local nonlinear Schr\"odinger and wave equations
	\citep{Dodson+Luhrmann+Mendelson_2020,Killip+Murphy+Visan_2019},
	the nonlocal Hartree interaction requires a substantially different analytic treatment.
	In the local power-type setting, the nonlinear terms can be expanded into
	pointwise products. By contrast, the convolution $|x|^{-4}\ast\cdot$ does not
	admit the same pointwise factorization. Its Fourier multiplier is comparable to $|\xi|^{-(d-4)}$, and hence the
	convolution corresponds to the regularizing operator $|\nabla|^{-(d-4)}$.
	Using Plancherel's theorem, we exploit these $d-4$ orders of regularization
	through a bilinear decomposition adapted to the regularity of the interacting
	terms. The resulting bilinear estimates replace the pointwise product
	estimates used in the local setting and combine Hardy--Littlewood--Sobolev,
	weighted radial Strichartz, and local smoothing estimates.
\end{remark}

We conclude the introduction by describing the main ideas of the proof. In the first part of the paper, we establish
the probabilistic estimates. Combining the weighted radial Strichartz estimate,
a local smoothing estimate that absorbs the derivative on the random evolution,
and Gaussian moment bounds, we show that for radial $f\in H^s$ with
$s>\frac{d}{2(d-1)}$, the forcing $F=e^{it\Delta}f^\omega$ almost surely lies
in the spacetime spaces used below. The condition $s>1-\frac{2}{d}$ provides the additional control required for
the energy increment.
In the second part of the paper we develop a stability theory for the forced
equation \eqref{forced1}, built on the deterministic scattering theory for
\eqref{HE1} \citep{Miao07, Miao11}. Treating $\mathcal{N}:=(|x|^{-4}\ast|v+F|^2)(v+F)-(|x|^{-4}\ast|v|^2)v$ as a
perturbation of the unforced equation and applying the deterministic global
spacetime bound, we show that a uniform bound $\sup_t\|v(t)\|_{H^1_x}\le E$ on the
maximal lifespan suffices for $v$ to be global and to scatter. This is the conditional criterion of
Theorem~\ref{thm:local}, which we verify unconditionally for $d\ge7$.
In the third part of the paper we derive this uniform energy bound through a double
bootstrap between the energy $\EE(T)=\sup_{t\in[0,T]}E[v(t)]$ and a Morawetz
spacetime norm $A(T)$. The control of $A(T)$ relies on two Morawetz-type estimates. The
first is a weighted Morawetz inequality, based on the convex weight
$a(x)=\langle x\rangle$ instead of the singular Lin--Strauss weight $|x|$,
which yields a weighted bound for $\nabla v$. The second is an interaction Morawetz inequality, which gives the global spacetime control
$\||\nabla|^{-\frac{d-3}{2}}|v|^2\|_{L^2_{t,x}}$, independent of the length of
the time interval. The nonlocal nonlinearity is coupled to both estimates
through the bilinear estimate of Lemma~\ref{lem:Hartree-reduction}.

This paper is organized as follows. Section~\ref{SECT:Prelim} presents the analytic and probabilistic estimates, including the Hartree bilinear estimate and the almost sure bounds for the randomized free evolution. Section~\ref{SECT:Forced} develops the local well-posedness and stability theory for the forced equation~\eqref{forced1} and proves Theorem~\ref{thm:local}. Section~\ref{Morawetz} establishes the weighted and interaction Morawetz estimates. Section~\ref{energy-bds} proves the uniform energy bound via a double-bootstrap argument and completes the proof of Theorem~\ref{thm:global}.

\section{Analytic and probabilistic estimates}\label{SECT:Prelim}
In this section, we will introduce notation and several analytic
tools that will be used throughout the paper. The deterministic
estimates are collected in Subsection~\ref{SS:det}, while the probabilistic estimates needed for the main results are presented in Subsection~\ref{SS:sto}.

For any time interval $I\subseteq\mathbb R$, we write
$L_t^qL_x^r(I\times\mathbb R^d)$ for $L^q(I;L^r(\mathbb R^d))$. More generally, for $a,q,r\in[1,\infty]$,
\[
\|f\|_{L_\omega^a L_t^q L_x^r}
:=\Bigl(\mathbb E\Bigl[\bigl\|\,\|f\|_{L_x^r(\mathbb R^d)}\bigr\|_{L_t^q}^{\,a}\Bigr]\Bigr)^{1/a}.
\] We denote by
$H_{\mathrm{rad}}^s(\mathbb R^d)$ the subspace of radial functions in
$H^s(\mathbb R^d)$.
A pair of exponents $(q,r)$ is called admissible if
$$
\frac{2}{q}+\frac{d}{r}=\frac d2,
\qquad
2\le q,r\le\infty,
\qquad
(q,r,d)\neq(2,\infty,2),
$$
and $(q',r')$ denotes the Hölder conjugate pair.
For later use, we introduce the following enlarged family of exponent pairs:
\begin{equation}\label{eq:Ad}
	\mathcal A_d:=
	\left\{
	(q,\rho):\
	2\le q<\infty,\ 
	\frac{2}{q}+\frac{d}{\rho}\le\frac d2
	\right\}.
\end{equation}
Equivalently, $(q,\rho)\in\mathcal A_d$ if and only if there exists an exponent
$r$ such that $(q,r)$ is admissible and $r\le \rho<\infty$.

We write $A\lesssim B$ whenever $A\le CB$ for some implicit constant
$C>0$ that may change from line to line, and $A\sim B$ whenever
$C_1B\le A\le C_2B$ for constants $C_1,C_2>0$ independent of the relevant
parameters.

\subsection{Deterministic results}\label{SS:det}
The Fourier transform on $\mathbb{R}^d$ is defined as
\begin{equation*}
	\hat f(\xi)\coloneqq(2\pi)^{-d/2}\int_{\mathbb{R}^d}
	f(x)\,e^{-ix\cdot\xi}\,\mathrm{d}x,
\end{equation*}
which allows us to define the fractional differentiation operators
$|\nabla|^s$ and $\langle\nabla\rangle^s$ through
\begin{equation*}
	\widehat{|\nabla|^{s}f}(\xi)\coloneqq|\xi|^s\hat f(\xi),
	\qquad
	\widehat{\langle\nabla\rangle^{s}f}(\xi)\coloneqq\langle\xi\rangle^s\hat f(\xi),
\end{equation*}
where $\langle\xi\rangle\coloneqq\sqrt{1+|\xi|^2}$. The homogeneous
and inhomogeneous Sobolev norms are then given by
\begin{align*}
	\|f\|_{\dot H_x^s(\mathbb{R}^d)}
	&\coloneqq
	\bigl\||\xi|^s\widehat f(\xi)\bigr\|_{L_\xi^2(\mathbb{R}^d)},\\
	\|f\|_{H_x^s(\mathbb{R}^d)}
	&\coloneqq
	\bigl\|\langle\xi\rangle^s\widehat f(\xi)\bigr\|_{L_\xi^2(\mathbb{R}^d)}.
\end{align*}
 
 Let $\varphi\in C_c^\infty(\mathbb{R}^d)$ be a radial bump function supported in
$\{|\xi|\le 2\}$ and satisfying $\varphi(\xi)=1$ for $|\xi|\le 1$. For each
dyadic number $N=2^k$, $k\in\mathbb{Z}$, define the Littlewood--Paley projection
operators by
\begin{align*}
	\widehat{P_{\le N}f}(\xi)
	&\coloneqq \varphi(\xi/N)\widehat f(\xi),\\
	\widehat{P_{>N}f}(\xi)
	&\coloneqq \bigl(1-\varphi(\xi/N)\bigr)\widehat f(\xi),\\
	\widehat{P_N f}(\xi)
	&\coloneqq \bigl(\varphi(\xi/N)-\varphi(2\xi/N)\bigr)\widehat f(\xi).
\end{align*}
Throughout, summation in $N$ is taken over
$N\in 2^{\mathbb{N}}=\{1,2,4,\ldots\}$, where $P_1$ denotes $P_{\le 1}$. Since $\varphi$ is radial, each $P_N$ preserves radial symmetry.
The Littlewood--Paley operators satisfy the following standard Bernstein estimates.

\begin{lem}[Bernstein estimates]\label{ber}
	Let $1\le p\le q\le \infty$ and $s\ge 0$. Then
	\begin{align*}
		\bigl\||\nabla|^s P_N u\bigr\|_{L_x^p(\mathbb{R}^d)}
		&\lesssim N^s \|P_N u\|_{L_x^p(\mathbb{R}^d)},\\
		\|P_N u\|_{L_x^q(\mathbb{R}^d)}
		&\lesssim N^{d(\frac{1}{p}-\frac{1}{q})}
		\|P_N u\|_{L_x^p(\mathbb{R}^d)}.
	\end{align*}
\end{lem}

We next recall two weighted estimates from
\citep{Killip+Murphy+Visan_2019}.
\begin{lem}\label{L.WY}
	\textup{(\citep[Lemma~2.2]{Killip+Murphy+Visan_2019})}
Let $1\le r\le m\le\infty$, $\beta \ge0$, and $\phi\in\mathcal{S}(\mathbb{R}^d)$. Then
\[
\bigl\|\langle x\rangle^{\beta}\bigl[|\phi|\ast|u|\bigr]\bigr\|_{L_x^m(\mathbb{R}^d)}
\lesssim
\bigl\|\langle x\rangle^{\beta} u\bigr\|_{L_x^r(\mathbb{R}^d)}.
\]
\end{lem}
\begin{lem}\label{lem:infinty}
	\textup{(\citep[Lemma~2.3]{Killip+Murphy+Visan_2019})}
	For $0\leq \beta \leq 1$ and $d<m<\infty$,
	\[
	\|\langle x\rangle^\beta u\|_{L_x^\infty(\mathbb{R}^d)}
	\lesssim
	\sum_{N\geq 1}
	N^{\frac{d}{m}}
	\|\langle x\rangle^\beta P_N u\|_{L_x^m(\mathbb{R}^d)} .
	\]
\end{lem}

We also collect several estimates for the Hartree nonlinearity. 

\begin{lem}\label{lem:Hardy}
		\textup{(Hardy's inequality \citep[Lemma A.2]{Tao2006})}
	Let $0\le s<\frac d2$. Then
	\[
	\bigl\||x|^{-s}f\bigr\|_{L_x^2(\mathbb R^d)}
	\lesssim
	\|f\|_{\dot H_x^s(\mathbb R^d)}.
	\]
\end{lem}
\begin{lem}\label{lem:1} 
	\textup{(Hardy-Littlewood-Sobolev inequality \citep[Theorem 4.3]{lieb01})}
	Let \( 0<\gamma<d \), \( 1 < p < q < \infty \), and 
$
	\frac{1}{p} + \frac{\gamma}{d} = 1 + \frac{1}{q}.
$
	Then, for any \( f \in L^p(\mathbb{R}^d) \), we have
	\[
\bigl\| |x|^{-\gamma} * f\bigr\|_{L_x^q(\mathbb{R}^d)} 	\lesssim \|f\|_{L_x^p(\mathbb{R}^d)}.
	\]
\end{lem}

We also recall the following identity for Riesz potentials.
\begin{lem}
	\label{lem:2.3}
	\textup{(\citep[Theorem~5.10]{lieb01})}
	Let $d\ge 5$ and $0<\alpha<\frac d2$. Then, for every
	$f\in C_c^\infty(\mathbb R^d)$,
\begin{equation}\label{eq:2.6}
	\bigl\||\nabla|^{-\alpha}\bigl(|f|^{2}\bigr)\bigr\|_{L^{2}(\mathbb{R}^{d})}^{2}
	=
	C_{d,\alpha}
	\iint_{\mathbb{R}^{d}\times\mathbb{R}^{d}}
	\frac{|f(x)|^{2}\,|f(y)|^{2}}{|x-y|^{d-2\alpha}}\,
	\mathrm{d}x\,\mathrm{d}y.
\end{equation}
\end{lem}

The following bilinear estimate distributes the regularizing effect of the
Hartree potential between its two factors. For scalar- or
$\mathbb C^d$-valued functions $g$ and $h$, define
\[
\mathfrak B_\gamma(g,h)(t)
:=
\int_{\mathbb R^d}
\bigl(|x|^{-\gamma}*g\bigr)(t,x)\cdot h(t,x)\,dx.
\]
\begin{lem}[Hartree bilinear estimate]\label{lem:Hartree-reduction}
	Let $d\ge1$, $0<\gamma<d$, and $1\le a,b,b'\le\infty$ satisfy
$
	\frac1b+\frac1{b'}=\frac1a.
$
Then, for every $0\le s\le d-\gamma$ and all sufficiently regular
spacetime functions $g$ and $h$ that are either both scalar-valued
or both $\mathbb C^d$-valued,
	\begin{equation}\label{eq:hartree-base}
		\bigl\|\mathfrak B_\gamma(g,h)\bigr\|_{L_t^a}
		\lesssim
		\bigl\||\nabla|^{-s}g\bigr\|_{L_t^bL_x^2}
		\bigl\||\nabla|^{-(d-\gamma-s)}h\bigr\|_{L_t^{b'}L_x^2}.
	\end{equation}
	
If, in addition, $0<\gamma\le d-1$, then, for every
$0\le s\le d-\gamma-1$ and all sufficiently regular scalar-valued
functions $g$ and $\mathbb C^d$-valued functions $h$,
\begin{equation}\label{eq:hartree-grad}
	\bigl\|\mathfrak B_\gamma(\nabla g,h)\bigr\|_{L_t^a}
	\lesssim
	\bigl\||\nabla|^{-s}g\bigr\|_{L_t^bL_x^2}
	\bigl\||\nabla|^{-(d-\gamma-1-s)}h\bigr\|_{L_t^{b'}L_x^2}.
\end{equation}
\end{lem}

\begin{proof}
	Since
	$
	\widehat{|x|^{-\gamma}}(\xi)
	=c_{d,\gamma}|\xi|^{-(d-\gamma)},
	$
	Parseval's identity gives, for each fixed $t$,
	\[
	\mathfrak B_\gamma(g,h)(t)
	=
	c_{d,\gamma}
	\int_{\mathbb R^d}
	\bigl(|\xi|^{-s}\widehat g(t,\xi)\bigr)
	\bigl(
	|\xi|^{-(d-\gamma-s)}\widehat h(t,-\xi)
	\bigr)\,d\xi.
	\]
	Thus, by the Cauchy--Schwarz inequality and Plancherel's identity,
	\[
	\bigl|\mathfrak B_\gamma(g,h)(t)\bigr|
	\lesssim
	\bigl\||\nabla|^{-s}g(t)\bigr\|_{L_x^2}
	\bigl\||\nabla|^{-(d-\gamma-s)}h(t)\bigr\|_{L_x^2}.
	\]
	Taking the $L_t^a$ norm and applying H\"older's inequality
	gives \eqref{eq:hartree-base}.
	
	For \eqref{eq:hartree-grad}, we have
	\[
	\mathfrak B_\gamma(\nabla g,h)(t)
	=
	\int_{\mathbb R^d}
	\nabla\bigl(|x|^{-\gamma}*g\bigr)(t,x)
	\cdot h(t,x)\,dx.
	\]
	Since $0\le s\le d-\gamma-1$, we apply
	\eqref{eq:hartree-base} to $\nabla g$ with exponent $s+1$. Then, by Plancherel's identity, we obtain
\[
\begin{aligned}
	\bigl\|\mathfrak B_\gamma(\nabla g,h)\bigr\|_{L_t^a}
	&\lesssim
	\bigl\||\nabla|^{-(s+1)}\nabla g\bigr\|_{L_t^bL_x^2}
	\bigl\||\nabla|^{-(d-\gamma-1-s)}h\bigr\|_{L_t^{b'}L_x^2}
	\\
	&=
	\bigl\||\nabla|^{-s}g\bigr\|_{L_t^bL_x^2}
	\bigl\||\nabla|^{-(d-\gamma-1-s)}h\bigr\|_{L_t^{b'}L_x^2}.
\end{aligned}
\]
	which proves \eqref{eq:hartree-grad}.
\end{proof}

We recall the following Strichartz estimates.
\begin{prop}\label{P:strichartz}
	\textup{(Strichartz estimates \citep[Theorem 1.2]{Keel98})}
For $d\ge2$ and admissible pairs $(q_j,r_j)$, $j=1,2$, the following
estimates hold on any time interval $I\subset\mathbb R$ with
$t_0\in\overline I$:
	\begin{align*}
		\bigl\|e^{it\Delta}f\bigr\|_{L_t^{q_1}L_x^{r_1}(\mathbb{R}\times\mathbb{R}^d)}
		&\lesssim \|f\|_{L_x^2(\mathbb{R}^d)},\\
		\biggl\|\int_{t_0}^{t}e^{i(t-s)\Delta}F(s)\,\mathrm{d}s\biggr\|_{L_t^{q_1}L_x^{r_1}(I\times\mathbb{R}^d)}
		&\lesssim \|F\|_{L_t^{q_2'}L_x^{r_2'}(I\times\mathbb{R}^d)}.
	\end{align*}
\end{prop}

In the radial setting, the linear propagator $e^{it\Delta}$ obeys the following bound, where
$P_{\mathrm{rad}}$ denotes the projection onto radial functions.
\begin{lem}\label{lem:rde}
	\textup{(Radial kernel estimate \citep[Lemma 2.8]{Killip08})}
	For every $t\neq 0$ and $x,y\in\mathbb{R}^d$ with $d\ge 2$,
	\begin{equation*}
		\bigl|e^{it\Delta}P_{\mathrm{rad}}(x,y)\bigr|
		\lesssim |t|^{-\frac12}\,|x|^{-\frac{d-1}{2}}\,|y|^{-\frac{d-1}{2}}.
	\end{equation*}
\end{lem}

Combining Lemma~\ref{lem:rde} with the standard dispersive estimate 
$\bigl|e^{it\Delta}(x,y)\bigr|\lesssim|t|^{-\frac{d}{2}}$ yields
\[
\bigl|e^{it\Delta}P_{\mathrm{rad}}(x,y)\bigr|
\lesssim |t|^{-\frac{2}{q}}\,|x|^{-\frac{dq-4}{2q}}\,|y|^{-\frac{dq-4}{2q}}
\qquad\text{for every } \tfrac{4}{d}\le q\le 4.
\]
A standard $TT^{*}$ argument then gives the following weighted
radial Strichartz estimate.
\begin{prop}\label{P:RS}
	\textup{(Weighted radial Strichartz estimate)}
	Let $d\ge 3$ and $2<q\le 4$. Then, for every $f\in L^2_{\mathrm{rad}}(\mathbb{R}^d)$,
	\[
	\bigl\||x|^{\frac{dq-4}{2q}}\,e^{it\Delta}f\bigr\|_{L_t^{q}L_x^{\infty}(\mathbb{R}\times\mathbb{R}^d)}
	\lesssim \|f\|_{L_x^2(\mathbb{R}^d)}.
	\]
\end{prop}
Interpolating Propositions~\ref{P:strichartz} and~\ref{P:RS} gives the
following weighted estimates.
\begin{lem}\label{C:WRS}
	Let $d\ge 3$, $2<q\le 4$, and $0\le\beta <\frac{dq-4}{2q}$. Then, for every $f\in L^2_{\mathrm{rad}}(\mathbb{R}^d)$,
	\[
	\bigl\||x|^{\beta}\,e^{it\Delta}f\bigr\|_{L_t^{q}L_x^{\frac{2dq}{dq-4-2q\beta}}(\mathbb{R}\times\mathbb{R}^d)}
	\lesssim \|f\|_{L_x^2(\mathbb{R}^d)}.
	\]
\end{lem}

The following local smoothing estimate
\citep{ConsSaut, Vega88} controls derivatives of the randomized linear evolution.
\begin{lem}[Local smoothing]\label{L:LS}
	For every $\varepsilon>0$,
	\[
	\bigl\|\langle x\rangle^{-\frac12-\varepsilon}\,e^{it\Delta}f\bigr\|_{L_{t,x}^{2}(\mathbb{R}\times\mathbb{R}^d)}
	\lesssim \|f\|_{\dot H_x^{-\frac12}(\mathbb{R}^d)}.
	\]
\end{lem}

\subsection{Almost sure bounds}\label{SS:sto}
The purpose of this subsection is to gather the probabilistic
spacetime estimates that will be used in the proof of
Theorem~\ref{thm:local} and \ref{thm:global}. Throughout, $f^{\omega}$ denotes the randomization
of $f$ introduced in Definition~\ref{def:random} on $\mathbb R^d$.  We first record the following moment bound of $f^{\omega}$.
\begin{lem}\textup{(Moment bound \citep[Lemma~2.8]{Killip+Murphy+Visan_2019})}\label{L:SF}
	For  $m\geq 1$,
	\[
	\mathbb{E}\bigl(|f^{\omega}|^{2m}\bigr)
	\lesssim
	\bigl(|\check{\psi}|\ast |f|^{2}\bigr)^{m}.
	\]
\end{lem}
Then, the following weighted spacetime estimate is the main tool of this subsection.

\begin{lem}\label{lem:gen-moment}
	Let $d\ge1$, $1\le q<\infty$, $2\le m<\infty$, $2\le r\le m$, and $\beta\ge0$.
	Then 
	\begin{equation}\label{eq:gen-estimate}
		\mathbb E\left[
		\|\langle x\rangle^\beta e^{it\Delta}f^\omega\|_{L_t^qL_x^m}^{q}
		\right]
	\lesssim
		\|\langle x\rangle^\beta e^{it\Delta}f\|_{L_t^qL_x^r}^{q}.
	\end{equation}
	
	Moreover, for any $p>2$, any $m\ge \frac{2p}{p-2}$, and any $2\le r_1,r_2\le m$, we have
	\begin{equation}\label{eq:gen-estimate-ball}
		\begin{aligned}
			\mathbb E\left[
			\|\langle x\rangle^\beta e^{it\Delta}f^\omega\|_{L_t^{\frac{2p}{p-2}}L_x^m}^{\frac{2p}{p-2}}
			\right]
		\lesssim
			\|e^{it\Delta}f\|_{L_t^{\frac{2p}{p-2}}L_x^{r_1}(\mathbb R\times B)}^{\frac{2p}{p-2}}
+
			\||x|^\beta e^{it\Delta}f\|_{L_t^{\frac{2p}{p-2}}L_x^{r_2}(\mathbb R\times B^c)}^{\frac{2p}{p-2}},
		\end{aligned}
	\end{equation}
	where $B:=\{x\in\mathbb R^d:|x|\le1\}$ and $B^c:=\mathbb R^d\setminus B$. All spacetime norms are taken over $\mathbb R\times\mathbb R^d$.
\end{lem}
\begin{proof}
	We first prove \eqref{eq:gen-estimate}. Fix $t\in\mathbb R$, applying Lemma~\ref{L:SF} to $e^{it\Delta}f$, we obtain
	\begin{align}
		\mathbb E
		\|\langle x\rangle^\beta e^{it\Delta}f^\omega\|_{L_x^m}^m
		&=
		\int_{\mathbb R^d}
		\langle x\rangle^{\beta m}
		\mathbb E\left[
		|e^{it\Delta}f^\omega(x)|^m
		\right]\,dx
		\nonumber
		\\
		&\lesssim
		\int_{\mathbb R^d}
		\langle x\rangle^{\beta m}
		\bigl(|\check{\psi}|*|e^{it\Delta}f|^2\bigr)^{m/2}(x)\,dx
		\nonumber
		\\
		&=
		\|\langle x\rangle^{2\beta}
		(|\check{\psi}|*|e^{it\Delta}f|^2)\|_{L_x^{m/2}}^{m/2}.
		\label{eq:weighted-moment-start}
	\end{align}
By Lemma~\ref{L.WY}, we get
	\begin{align}
		\|\langle x\rangle^{2\beta}
		(|\check{\psi}|*|e^{it\Delta}f|^2)\|_{L_x^{m/2}}
		\lesssim
		\|\langle x\rangle^{2\beta}|e^{it\Delta}f|^2\|_{L_x^{r/2}}
		=
		\|\langle x\rangle^\beta e^{it\Delta}f\|_{L_x^r}^2.
		\label{eq:weighted-convolution-use}
	\end{align}
	Combining \eqref{eq:weighted-moment-start} and \eqref{eq:weighted-convolution-use}, we obtain
	\begin{equation}\label{eq:m-pointwise-bound}
		\mathbb E
		\|\langle x\rangle^\beta e^{it\Delta}f^\omega\|_{L_x^m}^m
		\lesssim
		\|\langle x\rangle^\beta e^{it\Delta}f\|_{L_x^r}^m.
	\end{equation}

Next, for each fixed $t$, the function
$\langle x\rangle^\beta e^{it\Delta}f^\omega$
is an $L_x^m(\mathbb R^d)$-valued Gaussian random variable. Thus, by the
Kahane--Khintchine inequality for Gaussian random variables with values in Banach spaces
\citep{ledoux91}, we have, for every $1\le q<\infty$,
\begin{equation*}
	\left(
	\mathbb E
	\|\langle x\rangle^\beta e^{it\Delta}f^\omega\|_{L_x^m}^q
	\right)^{1/q}
	\lesssim
	\left(
	\mathbb E
	\|\langle x\rangle^\beta e^{it\Delta}f^\omega\|_{L_x^m}^m
	\right)^{1/m}.
\end{equation*}
Combining this with \eqref{eq:m-pointwise-bound}, we obtain
$$
	\mathbb E
	\|\langle x\rangle^\beta e^{it\Delta}f^\omega\|_{L_x^m}^q
	\lesssim
	\|\langle x\rangle^\beta e^{it\Delta}f\|_{L_x^r}^q.
$$
It follows that
	\begin{align*}
		\mathbb E\left[
		\|\langle x\rangle^\beta e^{it\Delta}f^\omega\|_{L_t^qL_x^m}^{q}
		\right]
		=
		\int_{\mathbb R}
		\mathbb E
		\|\langle x\rangle^\beta e^{it\Delta}f^\omega\|_{L_x^m}^{q}\,dt
		\lesssim
		\int_{\mathbb R}
		\|\langle x\rangle^\beta e^{it\Delta}f\|_{L_x^r}^{q}\,dt
		=
		C
		\|\langle x\rangle^\beta e^{it\Delta}f\|_{L_t^qL_x^r}^{q}.
	\end{align*}
	This proves \eqref{eq:gen-estimate}.

 Fix $t\in\mathbb R$. By Lemma~\ref{L:SF} and  elementary inequality
	$
	(a+b)^{m/2}\lesssim a^{m/2}+b^{m/2}, a,b\ge0,
	$
	we have
	\begin{align*}
		\mathbb E
		\|\langle x\rangle^\beta e^{it\Delta}f^\omega\|_{L_x^m}^m
		&\lesssim
		\|\langle x\rangle^{2\beta}
		(|\check{\psi}|*(\mathbf 1_B|e^{it\Delta}f|^2))\|_{L_x^{m/2}}^{m/2}
		\nonumber
		\\
		&\quad+
		\|\langle x\rangle^{2\beta}
		(|\check{\psi}|*(\mathbf 1_{B^c}|e^{it\Delta}f|^2))\|_{L_x^{m/2}}^{m/2}.
	\end{align*}
	Applying Lemma~\ref{L.WY} to the two terms yields
	\begin{equation}\label{eq:split-m-bound}
		\mathbb E
		\|\langle x\rangle^\beta e^{it\Delta}f^\omega\|_{L_x^m}^m
		\lesssim
		\|e^{it\Delta}f\|_{L_x^{r_1}(B)}^m
		+
		\||x|^\beta e^{it\Delta}f\|_{L_x^{r_2}(B^c)}^m,
	\end{equation}
	where we used $\langle x\rangle^\beta\lesssim1$ on $B$ and
	$\langle x\rangle^\beta\lesssim |x|^\beta$ on $B^c$.
	
	Since $m\ge \frac{2p}{p-2}$, \eqref{eq:split-m-bound} implies
$$
\mathbb E
\|\langle x\rangle^\beta e^{it\Delta}f^\omega\|_{L_x^m}^{\frac{2p}{p-2}}
\lesssim
\|e^{it\Delta}f\|_{L_x^{r_1}(B)}^{\frac{2p}{p-2}}
+
\||x|^\beta e^{it\Delta}f\|_{L_x^{r_2}(B^c)}^{\frac{2p}{p-2}}.
$$
	Thus, 
	\begin{align*}
		\mathbb E\left[
		\|\langle x\rangle^\beta e^{it\Delta}f^\omega\|_{L_t^{\frac{2p}{p-2}}L_x^m}^{\frac{2p}{p-2}}
		\right]
		&\lesssim
		\|e^{it\Delta}f\|_{L_t^{\frac{2p}{p-2}}L_x^{r_1}(\mathbb R\times B)}^{\frac{2p}{p-2}}+
		\||x|^\beta e^{it\Delta}f\|_{L_t^{\frac{2p}{p-2}}L_x^{r_2}(\mathbb R\times B^c)}^{\frac{2p}{p-2}}.
	\end{align*}
	This proves \eqref{eq:gen-estimate-ball}.
\end{proof}
\begin{remark}
	For $d=4$, \eqref{eq:gen-estimate} and \eqref{eq:gen-estimate-ball} reduce to
	\citep[Lemma~2.9]{Killip+Murphy+Visan_2019}, where the exponents satisfy $q\le m$
	because the $q$-th moment in $\omega$ is raised to the $m$-th by H\"older's inequality.
	We instead use the Kahane--Khintchine inequality, under which all moments of the
	$L^m_x$-valued Gaussian $\langle x\rangle^\beta e^{it\Delta}f^\omega$ are equivalent with
	a constant depending only on $q$ and $m$. Hence $q$ is arbitrary in $[1,\infty)$ and independent of $m$. The estimate also holds for all $d\ge1$, 
	which allows the probabilistic estimates from Lemma~\ref{lem:gen-moment} and Proposition~\ref{P:STE2} to be combined uniformly for the dimensions considered in
	Theorems~\ref{thm:local} and \ref{thm:global}.
\end{remark}

We now turn to the almost sure spacetime estimates.
\begin{prop}\label{P:STE2}
	Let $d\ge 5$ and $f\in L^2(\mathbb R^d)$, and let
	$\mathcal A\subset\mathcal A_d$ be a finite set, where $\mathcal A_d$ is
	defined in \eqref{eq:Ad}. Then
	\begin{equation}\label{E:STE2-finite-admissible}
		e^{it\Delta}f^\omega
		\in
		\bigcap_{(q,\rho)\in\mathcal A}
		L_t^qL_x^\rho \cap L_t^\infty L_x^2
	\end{equation}
	almost surely. Moreover, if $d=5$ and
	$f\in H^s(\mathbb R^5)$ for some $s>0$, then
	\begin{equation}\label{E:STE2-d5-L2Linf}
		e^{it\Delta}f^\omega\in L_t^2L_x^\infty
	\end{equation}
	almost surely. In addition, if $f\in H^s(\mathbb R^d)$ for some
	$s>\frac{d-2}{d}$, then
	\begin{equation}\label{E:STE2-Linfty-critical}
		e^{it\Delta}f^\omega
		\in
		L_t^\infty L_x^{\frac{2d}{d-2}}
	\end{equation}
	almost surely, where all spacetime norms are taken over
	$\mathbb R\times\mathbb R^d$.
\end{prop}
\begin{proof}
	Almost sure finiteness of the $L_t^\infty L_x^2$ norm follows from the
unitarity of the linear propagator on $L_x^2$ and the almost sure finiteness
of $\|f^\omega\|_{L_x^2}$:
$$
\|e^{it\Delta}f^\omega\|_{L_t^\infty L_x^2}
=
\|f^\omega\|_{L_x^2}
<\infty
\qquad\text{almost surely.}
$$
Next, let $(q,\rho)\in\mathcal A$. Since $\mathcal A\subset\mathcal A_d$,
there exists an exponent $r$ such that $(q,r)$ is admissible and
$r\le \rho<\infty$. By
	Lemma~\ref{lem:gen-moment} with $\beta=0$, followed by Proposition~\ref{P:strichartz}, we obtain
\begin{equation}\label{eq:random-strichartz-moment}
	\mathbb E
	\left(
	\|e^{it\Delta}f^\omega\|_{L_t^q L_x^\rho}^{q}
	\right)
	\lesssim
	\|e^{it\Delta}f\|_{L_t^q L_x^r}^{q}
	\lesssim
	\|f\|_{L_x^2}^{q}.
\end{equation}
	Thus
	$$
	\|e^{it\Delta}f^\omega\|_{L_t^q L_x^\rho}<\infty
	\qquad\text{almost surely.}
	$$
	Since $\mathcal A$ is finite, intersecting the corresponding
	probability-one events yields \eqref{E:STE2-finite-admissible}
	almost surely.
	
	We now prove the $L_t^2L_x^\infty$ assertion when $d=5$. Let
	$r=\frac{2d}{d-2}=\frac{10}{3}$, so that $(2,r)$ is admissible. Choose
	$m\ge r$ sufficiently large so that $\frac{d}{m}<s$. By Lemma~\ref{ber},
	$$
	\|e^{it\Delta}f^\omega\|_{L_t^2L_x^\infty}
	\lesssim
	\sum_{N\ge 1}
	N^{\frac{d}{m}}
	\|e^{it\Delta}P_N f^\omega\|_{L_t^2L_x^m}.
	$$
By	Lemma~\ref{lem:gen-moment} with $\beta=0$ and Proposition~\ref{P:strichartz},
	$$
	\begin{aligned}
		\left(
		\mathbb E
		\|e^{it\Delta}f^\omega\|_{L_t^2L_x^\infty}^{2}
		\right)^{1/2}
		&\lesssim
		\sum_{N\ge 1}
		N^{\frac{d}{m}}
		\left(
		\mathbb E
		\|e^{it\Delta}P_N f^\omega\|_{L_t^2L_x^m}^{2}
		\right)^{1/2}
		\\
		&\lesssim
		\sum_{N\ge 1}
		N^{\frac{d}{m}}
		\|e^{it\Delta}P_N f\|_{L_t^2L_x^r}
		\\
		&\lesssim
		\sum_{N\ge 1}
		N^{\frac{d}{m}}
		\|P_N f\|_{L_x^2}.
	\end{aligned}
	$$
	Since $\frac{d}{m}<s$, we have
	$$
	\sum_{N\ge 1}
	N^{\frac{d}{m}}
	\|P_N f\|_{L_x^2}
	=
	\sum_{N\ge 1}
	N^{\frac{d}{m}-s}
	N^s\|P_N f\|_{L_x^2}
	\lesssim
	\|f\|_{H_x^s}.
	$$
	Therefore
	$$
	\mathbb E
	\left(
	\|e^{it\Delta}f^\omega\|_{L_t^2L_x^\infty}^{2}
	\right)
	\lesssim
	\|f\|_{H_x^s}^{2},
	$$
	which gives \eqref{E:STE2-d5-L2Linf}
	almost surely.

	It remains to prove \eqref{E:STE2-Linfty-critical} almost surely.	
Fix a bounded interval $I\subset\mathbb{R}$ and let
$F:I\times\mathbb{R}^{d}\to\mathbb{C}$ be smooth enough. For every $t_{0}\in I$, 
\[
\|F\|_{L_{t}^{\infty}L_{x}^{\frac{2d}{d-2}}(I\times\mathbb{R}^{d})}
\le
\|F(t_{0})\|_{L_{x}^{\frac{2d}{d-2}}}
+\|\partial_{t}F\|_{L_{t}^{1}L_{x}^{\frac{2d}{d-2}}(I\times\mathbb{R}^{d})}.
\]
Averaging in $t_{0}\in I$ and applying H\"older in time reduces this to
\begin{equation}\label{E:FTC-step}
	\|F\|_{L_{t}^{\infty}L_{x}^{\frac{2d}{d-2}}(I\times\mathbb{R}^{d})}
	\le
	|I|^{\frac{2-d}{2d}}\,\|F\|_{L_{t,x}^{\frac{2d}{d-2}}(I\times\mathbb{R}^{d})}
	+|I|^{\frac{d+2}{2d}}\,\|\partial_{t}F\|_{L_{t,x}^{\frac{2d}{d-2}}(I\times\mathbb{R}^{d})}.
\end{equation}
Partitioning $\mathbb{R}$ into disjoint intervals of length $\delta>0$
and summing the $\tfrac{2d}{d-2}$-th power of \eqref{E:FTC-step} over
the partition gives
\begin{equation}\label{E:gen}
	\|F\|_{L_{t}^{\infty}L_{x}^{\frac{2d}{d-2}}(\mathbb{R}\times\mathbb{R}^{d})}^{\frac{2d}{d-2}}
	\;\lesssim\;
	\delta^{-1}\,\|F\|_{L_{t,x}^{\frac{2d}{d-2}}(\mathbb{R}\times\mathbb{R}^{d})}^{\frac{2d}{d-2}}
	+\delta^{\frac{d+2}{d-2}}\,\|\partial_{t}F\|_{L_{t,x}^{\frac{2d}{d-2}}(\mathbb{R}\times\mathbb{R}^{d})}^{\frac{2d}{d-2}}.
\end{equation}

Applying \eqref{E:gen} with $F=e^{it\Delta}P_N f^\omega$, and then using
Lemma~\ref{ber} together with \eqref{eq:random-strichartz-moment} applied to the pair
$\left(\frac{2d}{d-2},\frac{2d}{d-2}\right)\in\mathcal A_d$, we obtain
\begin{align*}
	\E\bigl( \| e^{it\Delta} P_{N} f^\omega \|_{L_t^\infty L_x^{\frac{2d}{d-2}}}^{\frac{2d}{d-2}}\bigr) & \lesssim  (\delta^{-1} + \delta^{\frac{d+2}{d-2}} N^{\frac{4d}{d-2}})\E\bigl( \|e^{it\Delta} P_{N} f^\omega\|_{L_{t,x}^{\frac{2d}{d-2}}}^{\frac{2d}{d-2}}\bigr) \\
	& \lesssim (\delta^{-1} + \delta^{\frac{d+2}{d-2}} N^{\frac{4d}{d-2}})\|P_Nf\|_{L^2_x}^{\frac{2d}{d-2}} \\
	& \lesssim (\delta^{-1}N^{-\frac{2ds}{d-2}}+ \delta^{\frac{d+2}{d-2}} N^{\frac{4d}{d-2}-\frac{2ds}{d-2}})\|f\|_{H_x^s}^{\frac{2d}{d-2}}.
\end{align*}
Optimizing in $\delta$, we get
\[
\E\bigl(\|e^{it\Delta}P_Nf^\omega\|_{L_t^\infty L_x^\frac{2d}{d-2}}^\frac{2d}{d-2}\bigr) \lesssim N^{2-\frac{2ds}{d-2}}\|f\|_{H_x^s}^\frac{2d}{d-2}.
\]
Thus
\[
\|e^{it\Delta} f^\omega\|_{L_\omega^{\frac{2d}{d-2}} L_t^\infty L_x^{\frac{2d}{d-2}}} \lesssim \sum_{N\geq 1} N^{\frac{d-2}{d}-s}\|f\|_{H_x^s} \lesssim \|f\|_{H_x^s}
\]
whenever $s>\frac{d-2}{d}$, which establishes \eqref{E:STE2-Linfty-critical} almost surely.
\end{proof}

\begin{prop}\label{P:STE}
	Let $d\geq5$, $s>\frac{d}{2(d-1)}$, and let
	$f\in H_{\mathrm{rad}}^s(\mathbb R^d)$. For $p$ sufficiently large, set
	\begin{equation}\label{def:M0}
		\mathcal{M}_0
		:=
		\left\{
		\frac{4d}{2d-j}:j\in\{4,5,6,8\}
		\right\}
		\cup
		\left\{
		\frac{2pd}{p(d-j)-4}:j\in\{2,3\}
		\right\}
		\cup
		\left\{
		\frac{2pd}{p(d-5)-4}:d>5
		\right\}.
	\end{equation}
	Then, for every $m\in\mathcal M_0$, one has
	\begin{equation}\label{E:nabla-common}
		\bigl\|
		\langle x\rangle^{\frac{3}{p}}
		\nabla e^{it\Delta}f^\omega
		\bigr\|_{L_t^{\frac{2p}{p-2}}L_x^m
			(\mathbb R\times\mathbb R^d)}
		<\infty
		\quad\text{almost surely}.
	\end{equation}
	Moreover, \eqref{E:nabla-common} remains valid with
	$m=\frac{2d}{d-5}$ when $d>5$, and with $m=\infty$ when $d=5$.
\end{prop}

\begin{proof}
	Fix a small parameter $\varepsilon>0$, to be chosen later, and define
	$$
	\theta
	:=
	\frac{p(1+2\varepsilon)+2}{p(d-1+2\varepsilon)}
	\in(0,1),
	\qquad
	q
	:=
	\frac{p(1+2\varepsilon)+2}
	{\left(\frac12+\varepsilon\right)(p-2)+3-d}.
	$$
	
	For $p$ sufficiently large and $\varepsilon>0$ sufficiently small,
we have $2<q\leq4$, and every finite exponent
$m\in\mathcal M_0$, as well as $m=\frac{2d}{d-5}$ when $d>5$,
satisfies
$
m\geq
\max\left\{
2,\frac{2}{1-\theta}
\right\}.
$
Hence, applying Lemma~\ref{lem:gen-moment} with
$r_1=2$ and $r_2=\frac{2}{1-\theta}$, we obtain
	\begin{equation}\label{E:split}
		\begin{aligned}
			&\mathbb E\Bigl[
			\bigl\|\langle x\rangle^{\frac{3}{p}}
			e^{it\Delta}\nabla P_Nf^\omega
			\bigr\|_{L_t^{\frac{2p}{p-2}}L_x^m}^{\frac{2p}{p-2}}
			\Bigr]
			\\
			&\qquad\lesssim
			\|e^{it\Delta}\nabla P_Nf\|_{L_t^{\frac{2p}{p-2}}L_x^{r_1}(\mathbb R\times B)}^{\frac{2p}{p-2}}
			+
			\bigl\||x|^{\frac{3}{p}}e^{it\Delta}\nabla P_Nf
			\bigr\|_{L_t^{\frac{2p}{p-2}}L_x^{r_2}(\mathbb R\times B^c)}^{\frac{2p}{p-2}}
				\\
		&\qquad:= I_N+II_N.
		\end{aligned}
	\end{equation} 
We first estimate $I_N$. By H\"older's inequality, Lemma~\ref{ber}, and Lemma~\ref{L:LS}, we have
	\begin{equation}\label{E:LS1}
		\begin{aligned}
		  I_N^{\frac{p-2}{2p}}
		&\lesssim
			\|e^{it\Delta}\nabla P_Nf\|_{L_{t,x}^2(\mathbb R\times B)}^{1-\frac{2}{p}}
			\|e^{it\Delta}\nabla P_Nf\|_{L_t^\infty L_x^2}^{\frac{2}{p}}
			\\
			&\lesssim
			\|\nabla P_Nf\|_{\dot H_x^{-\frac12}}^{1-\frac{2}{p}}
			\|\nabla P_Nf\|_{L_x^2}^{\frac{2}{p}}
			\\
			&\lesssim
			N^{\frac12+\frac1p}\|P_Nf\|_{L_x^2}.
		\end{aligned}
	\end{equation}
For $II_N$, Proposition~\ref{P:RS}, Lemma~\ref{L:LS}, and H\"older's inequality imply, for $\varepsilon>0$ sufficiently small and $p$ sufficiently large, that
	\begin{equation}\label{E:2nd-bd}
		\begin{aligned}
			 II_N^{\frac{p-2}{2p}}
			&\lesssim
			\bigl\||x|^{\frac{dq-4}{2q}}e^{it\Delta}\nabla P_Nf
			\bigr\|_{L_t^qL_x^\infty}^{\theta}
			\bigl\|\langle x\rangle^{-\frac12-\varepsilon}e^{it\Delta}\nabla P_Nf
			\bigr\|_{L_{t,x}^2}^{1-\theta}
			\\
			&\qquad\lesssim
			\|\nabla P_Nf\|_{L_x^2}^{\theta}
			\|\nabla P_Nf\|_{\dot H_x^{-\frac12}}^{1-\theta}
			\\
			&\qquad\lesssim
			N^{\frac12+\frac{\theta}{2}}\|P_Nf\|_{L_x^2}
			=
			N^{\frac{p(d+4\varepsilon)+2}{2p(d-1+2\varepsilon)}}
			\|P_Nf\|_{L_x^2}.
		\end{aligned}
	\end{equation}
Combining
	\eqref{E:split}, \eqref{E:LS1}, and \eqref{E:2nd-bd} yields
$$
\mathbb E\Bigl[
\bigl\|\langle x\rangle^{\frac{3}{p}}
e^{it\Delta}\nabla P_Nf^\omega
\bigr\|_{L_t^{\frac{2p}{p-2}}L_x^m}^{\frac{2p}{p-2}}
\Bigr]
\lesssim
\left(
N^{\frac{p(d+4\varepsilon)+2}{2p(d-1+2\varepsilon)}}
\|P_Nf\|_{L_x^2}
\right)^{\frac{2p}{p-2}}.
$$

Thus, we obtain
	$$
	\begin{aligned}
		\bigl\|\langle x\rangle^{\frac{3}{p}}
		\nabla e^{it\Delta}f^\omega
		\bigr\|_{L_\omega^{\frac{2p}{p-2}}L_t^{\frac{2p}{p-2}}L_x^m}
		&\lesssim
		\sum_{N\ge1}
		N^{\frac{p(d+4\varepsilon)+2}{2p(d-1+2\varepsilon)}}
		\|P_Nf\|_{L_x^2}
		&\lesssim
		\|f\|_{H_x^s},
	\end{aligned}
	$$
	provided that
	$$
	s>
	\frac{p(d+4\varepsilon)+2}{2p(d-1+2\varepsilon)}.
	$$
	Since $s>\frac{d}{2(d-1)}$, we choose
$\varepsilon>0$ sufficiently small and $p$
sufficiently large so that the above condition holds.
Consequently, \eqref{E:nabla-common} holds for every finite
$m\in\mathcal M_0$ and, when $d>5$, also for
$m=\frac{2d}{d-5}$.
	
	It remains to prove \eqref{E:nabla-common} for $m=\infty$ when $d=5$. Choose $M>d$ and $p$ sufficiently large so that
	$$
	M\ge \frac{2}{1-\theta}
	\qquad\text{and}\qquad
	s>
	\frac{p(d+4\varepsilon)+2}{2p(d-1+2\varepsilon)}
	+
	\frac{d}{M}.
	$$
	The preceding finite-exponent argument applies with $m=M$. Hence, by
	Lemma~\ref{lem:infinty},
	$$
	\begin{aligned}
		\bigl\|\langle x\rangle^{\frac{3}{p}}
		\nabla e^{it\Delta}f^\omega
		\bigr\|_{L_\omega^{\frac{2p}{p-2}}L_t^{\frac{2p}{p-2}}L_x^\infty}
		&\lesssim
		\sum_{N\ge1}
		N^{\frac{d}{M}}
		\bigl\|\langle x\rangle^{\frac{3}{p}}
		e^{it\Delta}\nabla P_Nf^\omega
		\bigr\|_{L_\omega^{\frac{2p}{p-2}}L_t^{\frac{2p}{p-2}}L_x^M}
		\\
		&\lesssim
		\sum_{N\ge1}
		N^{\frac{d}{M}}
		N^{\frac{p(d+4\varepsilon)+2}{2p(d-1+2\varepsilon)}}
		\|P_Nf\|_{L_x^2}
		\lesssim
		\|f\|_{H_x^s}.
	\end{aligned}
	$$
This completes the proof.
\end{proof}

\section{Well-posedness and scattering for the forced equation}\label{SECT:Forced}

In this section we study the Cauchy problem for the forced Hartree equation
in dimensions $d\geq 5$:
\begin{equation}\label{NLS-hartree-v1}
	\begin{cases}
		(i\partial_t+\Delta)v
		=
		\bigl(|x|^{-4}*|v+F|^2\bigr)(v+F),\\
		v(t_0,x)=v_0(x),
	\end{cases}
	\qquad (t,x)\in\mathbb{R}\times\mathbb{R}^d.
\end{equation}
where $t_0\in\mathbb{R}$ and the forcing term $F$ is a global solution to
the free Schrödinger equation
\begin{equation}\label{eq:F-free}
(i\partial_t+\Delta)F=0.
\end{equation}
We establish local well-posedness in the energy space $H^1(\mathbb R^d)$, 
together with a blow-up alternative and a scattering criterion.
We work with the spacetime norm
\begin{equation}\label{X-norm}
	\|v\|_{X(I)}
	:=
	\|v\|_{L_t^4L_x^{\frac{2d}{d-3}}(I\times\mathbb R^d)}
\end{equation}
and the forcing norm
\begin{equation}\label{Y-norm}
	\begin{aligned}
		\|F\|_{Y(I)}
		:={}&
		\|F\|_{L_t^{\frac{4p}{p+2}}L_x^{\frac{2d}{d-3}}(I\times\mathbb R^d)}
		+
		\|F\|_{L_t^{\frac{4p}{p+2}}L_x^{\frac{2dp}{(d-5)p+2}}(I\times\mathbb R^d)}
	+
		\|F\|_{L_t^{\frac{2p}{p-2}}W_x^{1,\frac{2d}{d-2}}(I\times\mathbb R^d)},
	\end{aligned}
\end{equation}
where $p$ is a fixed sufficiently large parameter.
\begin{prop}[Local well-posedness for the forced Hartree equation]\label{P:lwp_hartree}
	Let $d \geq 5$, $t_0 \in \mathbb{R}$, $v_0 \in H^1(\mathbb{R}^d)$ 
with $\|v_0\|_{H_x^1} \leq R$, and let $F$ solve \eqref{eq:F-free} 
on an open interval $I \ni t_0$ with $F(t_0) \in L_x^2(\mathbb{R}^d)$.
Assume that
\begin{equation}\label{lwp-condition-hartree}
	\|S(t-t_0)v_0\|_{X(I)} + \|F\|_{Y(I)} \leq \eta
\end{equation}
for some $\eta \leq \eta_0(R)$, where $\eta_0 = \eta_0(R) > 0$ is 
sufficiently small. Then the Cauchy problem \eqref{NLS-hartree-v1} 
admits a unique solution
\[
v \in C(I;\, H^1(\mathbb{R}^d)) \cap X(I).
\]
	Moreover, $v$ extends uniquely to a maximal-lifespan interval $I_{\max}\ni t_0$, and the following hold:
	\begin{enumerate}[label=\emph{(\roman*)}]
		\item\label{lwp:blowup}
		If $\sup I_{\max}<\infty$, then
		\[
		\|v\|_{X(t_0,\sup I_{\max})}=\infty.
		\]
		\item\label{lwp:scattering}
		If $\sup I_{\max}=\infty$ and $\|v\|_{X(t_0,\infty)}<\infty$, then $v$ scatters forward in time. More precisely, there exists $v_+\in H^1(\mathbb{R}^d)$ such that
		\[
		\|v(t)-S(t)v_+\|_{H_x^1}\to 0
		\qquad\text{as } t\to+\infty.
		\]
	\end{enumerate}
\end{prop}
\begin{proof}
Denote
$
\mathcal H(u):=\bigl(|x|^{-4}\ast |u|^2\bigr)u.
$ We define the solution map $\Gamma$ of \eqref{NLS-hartree-v1} by
\begin{equation}\label{Gamma-def}
	\Gamma(v)(t):=S(t-t_{0})v_{0}
	-i\int_{t_{0}}^{t}S(t-s)\,\mathcal{H}(v(s)+F(s))\,ds,
	\qquad t\in I.
\end{equation}
	Let $C_0 \geq 1$ be a fixed constant absorbing the implicit constants 
from Lemma~\ref{lem:1}, Proposition~\ref{P:strichartz}, and the Sobolev 
embedding. We 
define
\begin{equation}\label{ball}
	B_{\eta}(I):=\bigl\{v\in C_{t}H_{x}^{1}\cap X(I):\;
	\|v\|_{L_{t}^{\infty}H_{x}^{1}(I\times\mathbb{R}^{d})}\le 2C_{0}R,\;
	\|v\|_{X(I)}\le 2C_{0}\eta\bigr\},
\end{equation}
endowed with the metric induced by $L_{t}^{\infty}H_{x}^{1}\cap X(I)$.

We note the pairs $\left(\frac{4p}{3(p-2)},\frac{2dp}{p(d+1)+6}\right)$ and $\left(\frac{4p}{3p-2},\frac{2dp}{p(d+1)+2}\right)$ are dual admissible for $p>6$ and $p>2$, respectively. Using these pairs together with Lemma~\ref{lem:1} and 
H\"older's inequality, we obtain
\[
\begin{aligned}
	\bigl\|(|x|^{-4}\ast |F|^2)F\bigr\|_{L_t^1W_x^{1,2}}
	&\lesssim
	\bigl\||x|^{-4}\ast |F|^2\bigr\|_{L_t^{\frac{2p}{p+2}}L_x^d}
	\|F\|_{L_t^{\frac{2p}{p-2}}W_x^{1,\frac{2d}{d-2}}}
	\\
	&\quad+
	\bigl\||x|^{-4}\ast (F\nabla F)\bigr\|_
	{L_t^{\frac{4p}{3p-2}}L_x^{\frac{2d}{3}}}
	\|F\|_{L_t^{\frac{4p}{p+2}}L_x^{\frac{2d}{d-3}}}
	\\
	&\lesssim
	\|F\|^2_{L_t^{\frac{4p}{p+2}}L_x^{\frac{2d}{d-3}}}
	\|F\|_{L_t^{\frac{2p}{p-2}}W_x^{1,\frac{2d}{d-2}}}.
\end{aligned}
\]
Similarly, using the Sobolev embedding $W_x^{1,2}\hookrightarrow L_x^{\frac{2d}{d-2}}$,
\begin{align*}
	&\bigl\|(|x|^{-4}\ast |v|^2)v\bigr\|_
	{L_t^2W_x^{1,\frac{2d}{d+2}}}
	\\
	&\lesssim
	\bigl\||x|^{-4}\ast |v|^2\bigr\|_{L_t^2L_x^d}
	\|v\|_{L_t^\infty W_x^{1,2}}
	+
	\bigl\||x|^{-4}\ast (v\nabla v)\bigr\|_{L_t^2L_x^{\frac d2}}
	\|v\|_{L_t^\infty L_x^{\frac{2d}{d-2}}}
	\\
	&\lesssim
	\|v^2\|_{L_t^2L_x^{\frac d{d-3}}}
	\|v\|_{L_t^\infty W_x^{1,2}}
	+
	\|v\nabla v\|_{L_t^2L_x^{\frac d{d-2}}}
	\|v\|_{L_t^\infty W_x^{1,2}}
	\\
	&\lesssim
	\|v\|_{L_t^4L_x^{\frac{2d}{d-3}}}^{2}
	\|v\|_{L_t^\infty W_x^{1,2}}.
\end{align*}

By the fractional Leibniz rule, H\"older's inequality and Lemma~\ref{lem:1},
\begin{align}
	\bigl\|(|x|^{-4}\ast vF)F\bigr\|_
	{L_t^{\frac{4p}{3p-2}}W_x^{1,\frac{2dp}{p(d+1)+2}}}
	&\lesssim
	\bigl\||x|^{-4}\ast vF\bigr\|_{L_t^{\frac{4p}{p+2}}L_x^{\frac{2dp}{3p+2}}}
	\|F\|_{L_t^{\frac{2p}{p-2}}W_x^{1,\frac{2d}{d-2}}}\notag\\
	&\quad+
	\bigl\||x|^{-4}\ast(v\nabla F)\bigr\|_{L_t^{\frac{2p}{p-2}}L_x^{\frac d3}}
	\|F\|_{L_t^{\frac{4p}{p+2}}L_x^{\frac{2dp}{(d-5)p+2}}}\notag\\
	&\quad+
	\bigl\||x|^{-4}\ast(F\nabla v)\bigr\|_{L_t^{\frac{4p}{p+2}}L_x^{\frac{2dp}{3p+2}}}
	\|F\|_{L_t^{\frac{2p}{p-2}}L_x^{\frac{2d}{d-2}}}\notag\\
	&\lesssim
	\|vF\|_{L_t^{\frac{4p}{p+2}}L_x^{\frac{2dp}{p(2d-5)+2}}}
	\|F\|_{L_t^{\frac{2p}{p-2}}W_x^{1,\frac{2d}{d-2}}}\notag\\
	&\quad+
	\|v\nabla F\|_{L_t^{\frac{2p}{p-2}}L_x^{\frac d{d-1}}}
	\|F\|_{L_t^{\frac{4p}{p+2}}L_x^{\frac{2dp}{(d-5)p+2}}}\notag\\
	&\quad+
	\|F\nabla v\|_{L_t^{\frac{4p}{p+2}}L_x^{\frac{2dp}{p(2d-5)+2}}}
	\|F\|_{L_t^{\frac{2p}{p-2}}L_x^{\frac{2d}{d-2}}}.
	\label{E:vF-F-est}
\end{align}
The following term is treated in the same way:
\begin{align}
	\bigl\|(|x|^{-4}\ast|F|^2)v\bigr\|_
	{L_t^{\frac{4p}{3p-2}}W_x^{1,\frac{2dp}{p(d+1)+2}}}
	&\lesssim
	\bigl\||x|^{-4}\ast|F|^2\bigr\|_{L_t^{\frac{4p}{3p-2}}W_x^{1,\frac{2dp}{p+2}}}
	\|v\|_{L_t^\infty W_x^{1,2}}\notag\\
	&\lesssim
	\bigl\||F|^2\bigr\|_{L_t^{\frac{4p}{3p-2}}W_x^{1,\frac{2dp}{p(2d-7)+2}}}
	\|v\|_{L_t^\infty W_x^{1,2}}.\label{E:vF-F-est2}
\end{align}
Summing \eqref{E:vF-F-est} and \eqref{E:vF-F-est2} and applying H\"older's inequality and the
Sobolev embedding $W_x^{1,2}\hookrightarrow L_x^{\frac{2d}{d-2}}$, we obtain
\begin{equation*}
	\bigl\|(|x|^{-4}\ast vF)F\bigr\|_
	{L_t^{\frac{4p}{3p-2}}W_x^{1,\frac{2dp}{p(d+1)+2}}}
	+
	\bigl\|(|x|^{-4}\ast|F|^2)v\bigr\|_
	{L_t^{\frac{4p}{3p-2}}W_x^{1,\frac{2dp}{p(d+1)+2}}}
	\lesssim
	\|F\|_{L_t^{\frac{2p}{p-2}}W_x^{1,\frac{2d}{d-2}}}
	\|F\|_{L_t^{\frac{4p}{p+2}}L_x^{\frac{2dp}{(d-5)p+2}}}
	\|v\|_{L_t^\infty W_x^{1,2}}.
\end{equation*}
The remaining two terms are estimated similarly:
\begin{align*}
	\bigl\|(|x|^{-4}\ast|v|^2)F\bigr\|_
	{L_t^{\frac{4p}{3(p-2)}}W_x^{1,\frac{2dp}{p(d+1)+6}}}
	&\lesssim
	\bigl\||x|^{-4}\ast|v|^2\bigr\|_{L_t^{\frac{4p}{p-2}}L_x^{\frac{2dp}{3(p+2)}}}
	\|F\|_{L_t^{\frac{2p}{p-2}}W_x^{1,\frac{2d}{d-2}}}\\
	&\quad+
	\bigl\||x|^{-4}\ast(v\nabla v)\bigr\|_{L_t^{\frac{4p}{p-2}}L_x^{\frac{2dp}{5p+6}}}
	\|F\|_{L_t^{\frac{2p}{p-2}}L_x^{\frac{2d}{d-4}}}\\
	&\lesssim
	\|v^2\|_{L_t^{\frac{4p}{p-2}}L_x^{\frac{2dp}{p(2d-5)+6}}}
	\|F\|_{L_t^{\frac{2p}{p-2}}W_x^{1,\frac{2d}{d-2}}}\\
	&\quad+
	\|v\nabla v\|_{L_t^{\frac{4p}{p-2}}L_x^{\frac{2dp}{p(2d-3)+6}}}
	\|F\|_{L_t^{\frac{2p}{p-2}}L_x^{\frac{2d}{d-4}}},
\end{align*}
and, proceeding as before,
\begin{align*}
	\bigl\|(|x|^{-4}\ast vF)v\bigr\|_
	{L_t^{\frac{4p}{3(p-2)}}W_x^{1,\frac{2dp}{p(d+1)+6}}}
	&\lesssim
	\bigl\||x|^{-4}\ast(vF)\bigr\|_{L_t^{\frac{4p}{3(p-2)}}L_x^{\frac{2dp}{p+6}}}
	\|v\|_{L_t^\infty W_x^{1,2}}\\
	&\quad+
	\bigl\||x|^{-4}\ast(v\nabla F)\bigr\|_{L_t^{\frac{4p}{3(p-2)}}L_x^{\frac{2dp}{3(p+2)}}}
	\|v\|_{L_t^\infty L_x^{\frac{2d}{d-2}}}\\
	&\quad+
	\bigl\||x|^{-4}\ast(F\nabla v)\bigr\|_{L_t^{\frac{2p}{p-2}}L_x^{\frac d2}}
	\|v\|_{L_t^\infty L_x^2}^{\frac2p}
	\|v\|_{L_t^4L_x^{\frac{2d}{d-3}}}^{1-\frac2p}\\
	&\lesssim
	\|vF\|_{L_t^{\frac{4p}{3(p-2)}}L_x^{\frac{2dp}{p(2d-7)+6}}}
	\|v\|_{L_t^\infty W_x^{1,2}}\\
	&\quad+
	\|v\nabla F\|_{L_t^{\frac{4p}{3(p-2)}}L_x^{\frac{2dp}{p(2d-5)+6}}}
	\|v\|_{L_t^\infty W_x^{1,2}}\\
	&\quad+
	\|F\nabla v\|_{L_t^{\frac{2p}{p-2}}L_x^{\frac d{d-2}}}
	\|v\|_{L_t^\infty L_x^2}^{\frac2p}
	\|v\|_{L_t^4L_x^{\frac{2d}{d-3}}}^{1-\frac2p}.
\end{align*}
Summing and using \eqref{X-norm}, we conclude
\begin{equation*}
	\bigl\|(|x|^{-4}\ast|v|^2)F\bigr\|_
	{L_t^{\frac{4p}{3(p-2)}}W_x^{1,\frac{2dp}{p(d+1)+6}}}
	+
	\bigl\|(|x|^{-4}\ast vF)v\bigr\|_
	{L_t^{\frac{4p}{3(p-2)}}W_x^{1,\frac{2dp}{p(d+1)+6}}}
	\lesssim
	\|F\|_{L_t^{\frac{2p}{p-2}}W_x^{1,\frac{2d}{d-2}}}
	\|v\|_{L_t^\infty L_x^2}^{\frac2p}
	\|v\|_X^{1-\frac2p}
	\|v\|_{L_t^\infty W_x^{1,2}}.
\end{equation*}

For $v \in B_\eta(I)$, applying Proposition~\ref{P:strichartz} and 
the Sobolev embedding 
$W^{1,\frac{2d}{d-1}} \hookrightarrow L^{\frac{2d}{d-3}}$,
together with the estimates above, we obtain
\begin{equation}\label{Gamma-Linf-H1}
\|\Gamma(v)\|_{L_{t}^{\infty}H_{x}^{1}}\lesssim R+\eta^{3}+\eta^{2} R+R^{1+\frac{2}{p}}\eta^{2-\frac{2}{p}}\le 2C_{0}R,
\end{equation}
and
\begin{equation}\label{Gamma-X}
	\|\Gamma(v)\|_{X(I)}\lesssim \eta+\eta^{3}+\eta^{2} R+R^{1+\frac{2}{p}}\eta^{2-\frac{2}{p}}\le 2C_{0}\eta,
\end{equation}
	provided $\eta \leq \eta_0(R)$ is chosen sufficiently small. 
Hence $\Gamma(B_\eta(I)) \subseteq B_\eta(I)$.

For $v_1,v_2\in B_\eta(I)$, set $w:=v_1-v_2$. Using the identity
\[
\begin{aligned}
	\mathcal H(v_1+F)-\mathcal H(v_2+F)
	&=\bigl(|x|^{-4}*|v_1+F|^2\bigr)w\\
	&\quad+\frac12
	\bigl[|x|^{-4}*
	\bigl((v_1+v_2+2F)\overline w\bigr)\bigr](v_2+F)\\
	&\quad+\frac12
	\bigl[|x|^{-4}*
	\bigl(\overline{(v_1+v_2+2F)}\,w\bigr)\bigr](v_2+F),
\end{aligned}
\]
together with the estimates above, we obtain
\[
\|\Gamma(v_1)-\Gamma(v_2)\|_{L_t^\infty H_x^1\cap X(I)}
\le
C_0\left(
\eta^2+R\eta+R^{\frac2p}\eta^{2-\frac2p}
\right)
\|v_1-v_2\|_{L_t^\infty H_x^1\cap X(I)}.
\]
Choosing $\eta_0(R)>0$ sufficiently small so that
\[
C_0\left(
\eta_0(R)^2
+R\eta_0(R)
+R^{\frac2p}\eta_0(R)^{2-\frac2p}
\right)
\le\frac12,
\]
we conclude that $\Gamma$ is a strict contraction on $B_\eta(I)$
whenever $\eta\le\eta_0(R)$. The Banach fixed-point theorem therefore
yields a unique $v\in B_\eta(I)$ satisfying $\Gamma(v)=v$, which solves
\eqref{NLS-hartree-v1} on $I$.
	
A standard continuation argument then yields a unique maximal-lifespan
solution on an open interval $I_{\max}\ni t_0$.
We now prove \ref{lwp:blowup}. Suppose toward a contradiction that
$\sup I_{\max}<\infty$ and
\[
\|v\|_{X((t_0,\sup I_{\max})\times\mathbb{R}^d)}<\infty.
\]
Note that since $F$ solves the free Schr\"odinger equation, $u:=v+F$
solves the unforced Hartree equation $(i\partial_t+\Delta)u=\mathcal H(u)$
on $I_{\max}\times\mathbb{R}^d$ with $u(t_0)=v_0+F(t_0)$. Thus, by mass
conservation and the triangle inequality,
\[
\|v\|_{L_t^\infty L_x^2((t_0,\sup I_{\max})\times\mathbb{R}^d)}
\lesssim \|v_0\|_{L_x^2}+\|F(t_0)\|_{L_x^2}.
\]
Decomposing $(t_0,\sup I_{\max})$ into finitely many subintervals $J$ on
which $\|v\|_{X(J)}+\|F\|_{Y(J)}$ is sufficiently small and iterating the
nonlinear estimates above, we obtain
\[
\|v\|_{L_t^\infty H_x^1((t_0,\sup I_{\max})\times\mathbb{R}^d)}<\infty.
\]
Fix $\tau \in (t_0, \sup I_{\max})$. By the Duhamel formula with initial
time $\tau$,
\[
v(t)=S(t-\tau)v(\tau)
-i\int_\tau^t S(t-s)\mathcal H(v+F)(s)\,ds,
\qquad t\in[\tau,\sup I_{\max}).
\]
Applying  Proposition~\ref{P:strichartz} together with the nonlinear estimates established above,
and using
\[
\|v\|_{X([\tau,\sup I_{\max}))}
+\|F\|_{Y([\tau,\sup I_{\max}))}
\longrightarrow 0
\qquad\text{as }\tau\to\sup I_{\max},
\]
we obtain
\[
\bigl\|S(\cdot-\tau)v(\tau)\bigr\|_{X([\tau,\sup I_{\max}))}
+\|F\|_{Y([\tau,\sup I_{\max}))}
\longrightarrow 0
\qquad\text{as }\tau\to\sup I_{\max}.
\]
In particular, for $\tau$ sufficiently close to $\sup I_{\max}$ the
smallness condition
\[
\bigl\|S(\cdot-\tau)v(\tau)\bigr\|_{X([\tau,\sup I_{\max}))}
+\|F\|_{Y([\tau,\sup I_{\max}))}
\le \eta_0\bigl(\|v(\tau)\|_{H_x^1}\bigr)
\]
is satisfied. The local theory then extends $v$ beyond 
$\sup I_{\max}$, contradicting the maximality of $I_{\max}$. This proves 
\ref{lwp:blowup}.

Finally, suppose $\sup I_{\max}=\infty$ and
$\|v\|_{X(t_0,\infty)}<\infty$. The same decomposition argument as above yields
$v\in L_t^\infty H_x^1(t_0,\infty)$. For $t>s>t_0$, the Duhamel formula and
the group property of $S(t)$ give
\[
S(-t)v(t)-S(-s)v(s)
=
-i\int_s^t S(-\sigma)H(v+F)(\sigma)\,d\sigma.
\]
Applying  Proposition~\ref{P:strichartz}  together with the nonlinear estimates established above
and the absolute continuity of the relevant spacetime norms, the Duhamel term on
the right-hand side tends to zero in $H_x^1$ as $s,t\to\infty$.
Hence $\{S(-t)v(t)\}_{t\ge t_0}$ is Cauchy in $H_x^1(\mathbb{R}^d)$.
Therefore, there exists $v_+\in H^1(\mathbb{R}^d)$ such that
\[
S(-t)v(t)\longrightarrow v_+
\quad\text{in }H_x^1
\]
as $t\to+\infty$. Since $S(t)$ is unitary on $H_x^1$,
\[
\|v(t)-S(t)v_+\|_{H_x^1}
=
\|S(-t)v(t)-v_+\|_{H_x^1}
\longrightarrow 0
\qquad\text{as }t\to+\infty.
\]
This proves \ref{lwp:scattering}.
	\end{proof}

We next establish a stability theory for the forced Hartree equation, beginning with a short-time stability result.

\begin{lem}[Short-time stability for the forced Hartree equation]
	\label{L:short-perturb-hartree}
	Let $d\ge 5$, let $I\ni t_0$ be a time interval, and let
	$v_0\in H^1(\mathbb R^d)$ satisfy
$
	\|v_0\|_{H_x^1}\le E.
$
	Suppose that
	$u\in C(I;H^1(\mathbb R^d))\cap X(I)$ solves
	\begin{equation}\label{eq:unforced-hartree}
		(i\partial_t+\Delta)u=\mathcal H(u)
	\end{equation}
	on $I\times\mathbb R^d$ and satisfies
	\[
	\sup_{t\in I}\|u(t)\|_{H_x^1}\le E,
	\qquad
	\|v_0-u(t_0)\|_{H_x^1}\le\varepsilon,
	\qquad
	\|u\|_{X(I)}\le\delta
	\]
	for some $0<\varepsilon\le\varepsilon_0$. Suppose also that $F$ solves
\eqref{eq:F-free}
	on $I$, with $F(t_0)\in L_x^2(\mathbb R^d)$, and satisfies
	\[
	\|F\|_{Y(I)}\le\varepsilon.
	\]
	
	Then, for $\varepsilon_0$ and $\delta$ sufficiently small depending on
	$E$, the maximal-lifespan solution $v$ to
	\begin{equation}\label{eq:forced-hartree}
		(i\partial_t+\Delta)v=\mathcal H(v+F),
		\qquad
		v(t_0)=v_0,
	\end{equation}
	with maximal-lifespan interval $I_{\max}$ satisfies
	$
	I\subset I_{\max}
$
	and
\begin{equation}\label{eq:short-time-stability}
	\|v-u\|_{L_t^\infty H_x^1\cap X(I)}
	\le C(E)\varepsilon.
\end{equation}
\end{lem}

\begin{proof}
	By time translation, we  assume that $t_0=0$. It suffices to consider
	the forward-time interval, since the backward-time argument is identical.
	All spacetime norms below are taken over an arbitrary compact subinterval
	of $I\cap I_{\max}$ containing $0$.
	
	Set $w:=v-u$. Then
	\[
	(i\partial_t+\Delta)w
	=
	\mathcal H(u+w+F)-\mathcal H(u),
	\qquad
	w(0)=v_0-u(0).
	\]
	By the Duhamel formula, Proposition~\ref{P:strichartz}, and the nonlinear
	estimates established in the proof of
	Proposition~\ref{P:lwp_hartree}, we obtain
	\[
	\begin{aligned}
		\|w\|_{L_t^\infty H_x^1\cap X}
		\lesssim{}&
		\varepsilon
		+\bigl(\delta^2+\|w\|_{L_t^\infty H_x^1\cap X}^{2}\bigr)
		\|w\|_{L_t^\infty H_x^1\cap X}\\
		&+\varepsilon
		\bigl(E+\|w\|_{L_t^\infty H_x^1\cap X}\bigr)^{1+\frac2p}
		\bigl(\delta+\|w\|_{L_t^\infty H_x^1\cap X}\bigr)^{1-\frac2p}\\
		&+\varepsilon^2
		\bigl(E+\|w\|_{L_t^\infty H_x^1\cap X}\bigr)
		+\varepsilon^3.
	\end{aligned}
	\]
	Since $p>2$, we  first choose $\delta>0$ sufficiently small depending
	on $E$ and then choose $\varepsilon_0>0$ sufficiently small depending
	on $E$. A standard continuity argument then yields
	\[
	\|w\|_{L_t^\infty H_x^1\cap X}
	\le C(E)\varepsilon
	\]
	uniformly over all compact subintervals of
	$I\cap I_{\max}$ containing $0$. Consequently,
	\[
	\|v\|_{X}
	\le
	\|u\|_{X(I)}+\|w\|_{X}
	\le
	\delta+C(E)\varepsilon.
	\]
The blow-up alternative in Proposition~\ref{P:lwp_hartree}, together
with the analogous backward-time argument, therefore yields
$
I\subset I_{\max}
$
and \eqref{eq:short-time-stability}.
\end{proof}

Next we derive a long-time stability result by iterating  Lemma \ref{L:short-perturb-hartree} over a finite partition of the interval.

\begin{lem}[Long-time stability for the forced Hartree equation]
	\label{L:long-stability-hartree}
	Let $d\ge 5$, let $I\ni t_0$ be a time interval, and let
	$v_0\in H^1(\mathbb R^d)$ with $\|v_0\|_{H_x^1}\le E$. Let $u$ be the
	solution to the unforced Hartree equation \eqref{eq:unforced-hartree}
	with $u(t_0)=v_0$, and suppose
	\[
	\|u\|_{X(I)}\le L.
	\]
	
	Then there exists $\varepsilon_1=\varepsilon_1(E,L)>0$ such that if
	$F$ solves \eqref{eq:F-free} on $I$, with
	$F(t_0)\in L_x^2(\mathbb R^d)$, and satisfies
	\begin{equation}\label{eq:F-smallness-long}
		\|F\|_{Y(I)}\le\varepsilon
	\end{equation}
	for some $0<\varepsilon\le\varepsilon_1$, then there exists a unique
	solution $v$ to \eqref{eq:forced-hartree} with $v(t_0)=v_0$ on
	$I\times\mathbb R^d$. Moreover,
	\begin{equation}\label{eq:long-time-stability}
		\|v-u\|_{L_t^\infty H_x^1\cap X(I)}
		\le C(E,L)\varepsilon,
		\qquad
		\|v\|_{X(I)}
		\le C(E,L)<\infty.
	\end{equation}
\end{lem}

\begin{proof}
	By the local theory, it suffices to establish the $X(I)$-bound for $v$
	as an a priori estimate. By conservation of mass and energy, we have
$
	\|u\|_{L_t^\infty H_x^1(I\times\mathbb R^d)}
	\le C_0(E).
$
	
 Choose $\delta=\delta(2C_0(E))$ as in
	Lemma~\ref{L:short-perturb-hartree} and divide the forward-time part of
	$I$ into $J=J(E,L)$ subintervals $I_j=[t_j,t_{j+1}]$ such that
	\[
	\|u\|_{X(I_j)}\le\delta
	\qquad\text{for each }j.
	\]
	
	Define
	$
	C_0=1,
$ $
	C_{j+1}=C(E)(C_j+1).
$
	Since $J$ depends only on $E$ and $L$, we may choose
	$\varepsilon_1=\varepsilon_1(E,L)>0$ sufficiently small so that
$
	\sup_{0\le j\le J}(C_j+1)\varepsilon_1
	\le\varepsilon_0(2C_0(E))
$
	and
$
	\sup_{0\le j\le J}C_j\varepsilon_1
	\le C_0(E).
$
	
	We claim that
	\begin{equation}\label{eq:inductive-endpoint-bound}
		\|v(t_j)-u(t_j)\|_{H_x^1}\le C_j\varepsilon,
		\qquad
		\|v(t_j)\|_{H_x^1}\le2C_0(E)
	\end{equation}
	for each $j$. The claim holds for $j=0$ since
	$v(t_0)=u(t_0)=v_0$. Suppose that it holds at $j$. Since
	\[
	\|v(t_j)-u(t_j)\|_{H_x^1}
	\le C_j\varepsilon
	\le(C_j+1)\varepsilon
	\]
	and
	\[
	\|F\|_{Y(I_j)}
	\le\varepsilon
	\le(C_j+1)\varepsilon,
	\]
	Lemma~\ref{L:short-perturb-hartree}, applied on $I_j$ with perturbation
	size $(C_j+1)\varepsilon$, gives
	\[
	\|v-u\|_{L_t^\infty H_x^1\cap X(I_j)}
	\le C(E)(C_j+1)\varepsilon
	=C_{j+1}\varepsilon.
	\]
	Hence
	\[
	\|v(t_{j+1})-u(t_{j+1})\|_{H_x^1}
	\le C_{j+1}\varepsilon.
	\]
	Moreover,
	\[
	\|v(t_{j+1})\|_{H_x^1}
	\le C_0(E)+C_{j+1}\varepsilon
	\le2C_0(E).
	\]
	
	Combining the estimates on the finitely many intervals $I_j$, we obtain
	\[
	\|v-u\|_{L_t^\infty H_x^1\cap X}
	\le C(E,L)\varepsilon
	\]
	on the forward-time part of $I$. The same argument backward in time
	gives
	\[
	\|v-u\|_{L_t^\infty H_x^1\cap X(I)}
	\le C(E,L)\varepsilon.
	\]
	Consequently,
	\[
	\|v\|_{X(I)}
	\le
	\|u\|_{X(I)}+\|v-u\|_{X(I)}
	\le C(E,L)<\infty.
	\]
	This proves \eqref{eq:long-time-stability} and completes the proof.
\end{proof}
We now combine Lemma~\ref{L:long-stability-hartree} with
Proposition~\ref{P:lwp_hartree} to obtain a scattering criterion for
the forced Hartree equation under an a priori $H_x^1$ bound.
\begin{prop}[$H^1$ bounds imply scattering for the forced Hartree equation]
	\label{P:H1-implies-scattering-hartree}
	Let $d\ge 5$, let $t_0\in\mathbb R$, and let $F$ solve
	$
	(i\partial_t+\Delta)F=0
	$
	on $[t_0,\infty)$ with $F(t_0)\in L_x^2(\mathbb R^d)$ and
	$
	\|F\|_{Y((t_0,\infty))}<\infty.
	$
	Let
	$
	v:I_{\max}\times\mathbb R^d\to\mathbb C
	$
	be the maximal-lifespan solution to
	\[
	(i\partial_t+\Delta)v=\mathcal H(v+F),
	\qquad
	v(t_0)=v_0\in H^1(\mathbb R^d).
	\]
	Suppose that
	\begin{equation}\label{eq:uniform-H1-bound-forced}
		\sup_{t\in I_{\max}\cap[t_0,\infty)}
		\|v(t)\|_{H_x^1}\le E<\infty .
	\end{equation}
	Then $\sup I_{\max}=\infty$ and $v$ scatters forward in time. More precisely,
	there exists $v_+\in H^1(\mathbb R^d)$ such that
\[
\|v(t)-S(t)v_+\|_{H_x^1}\to0
\qquad\text{as }t\to+\infty.
\]
\end{prop}

\begin{proof}
	By Proposition~\ref{P:lwp_hartree}, it suffices to show that
	\begin{equation}\label{eq:global-X-bound-goal}
		\|v\|_{X((t_0,\sup I_{\max}))}<\infty .
	\end{equation}
	Indeed, the blow-up alternative then gives $\sup I_{\max}=\infty$, and the
	scattering criterion in Proposition~\ref{P:lwp_hartree} yields forward
	scattering.
	
	We prove \eqref{eq:global-X-bound-goal}. By the global theory for the
	unforced Hartree equation in $H^1$, there exists $L(E)<\infty$ such that,
	for every solution $u$ to
	$
	(i\partial_t+\Delta)u=\mathcal H(u)
	$
	with $\|u(s)\|_{H_x^1}\le E$ at some time $s$, one has
	\begin{equation}\label{eq:unforced-global-X-bound}
		\|u\|_{X(\mathbb R)}\le L(E).
	\end{equation}
	
	Let $\varepsilon_1=\varepsilon_1(E,L(E))>0$ be as in
	Lemma~\ref{L:long-stability-hartree}. Since
	$\|F\|_{Y((t_0,\infty))}<\infty$, we divide
	$[t_0,\sup I_{\max})$ into finitely many consecutive intervals
	\[
	[t_0,\sup I_{\max})=\bigcup_{j=0}^{J-1} I_j,
	\qquad
	I_j=[t_j,t_{j+1}),
	\]
	such that
	\begin{equation}\label{eq:F-small-on-subintervals}
		\|F\|_{Y(I_j)}\le \varepsilon_1
		\qquad\text{for each } j .
	\end{equation}
	Here $J<\infty$ depends only on $\varepsilon_1$ and
	$\|F\|_{Y((t_0,\infty))}$.
	
	Fix $j$. Let $u_j$ be the solution to the unforced Hartree equation with
	initial data
	$
	u_j(t_j)=v(t_j).
	$
	By \eqref{eq:uniform-H1-bound-forced},
	\[
	\|u_j(t_j)\|_{H_x^1}
	=
	\|v(t_j)\|_{H_x^1}
	\le E.
	\]
	Thus, by \eqref{eq:unforced-global-X-bound},
	\[
	\|u_j\|_{X(I_j)}\le \|u_j\|_{X(\mathbb R)}\le L(E).
	\]
	Since $v$ and $u_j$ have the same initial data at time $t_j$, and since
	\eqref{eq:F-small-on-subintervals} holds, Lemma~\ref{L:long-stability-hartree}
	applies on $I_j$. Consequently,
	\begin{equation}\label{eq:X-bound-each-piece}
		\|v\|_{X(I_j)}\le C(E,L(E))=:C(E)
		\qquad\text{for each } j .
	\end{equation}
	
	Summing over the finitely many intervals gives
	\[
	\begin{aligned}
		\|v\|_{X((t_0,\sup I_{\max}))}^4
		&=
		\sum_{j=0}^{J-1}\|v\|_{X(I_j)}^4
		&\le
		JC(E)^4
		<\infty .
	\end{aligned}
	\]
	This proves \eqref{eq:global-X-bound-goal}.
\end{proof}

Combining Propositions~\ref{P:lwp_hartree}
and~\ref{P:H1-implies-scattering-hartree}, we now prove
Theorem~\ref{thm:local}.

\begin{proof}[Proof of Theorem~\ref{thm:local}]
	Let $d\ge5$ and let
	$f\in H_{\mathrm{rad}}^s(\mathbb R^d)$ with
$
	\frac{d}{2(d-1)}<s<1.
$
	Define
	\begin{equation}\label{def:AY}
		\mathcal A_Y
		:=
		\left\{
		\left(\frac{4p}{p+2},\frac{2d}{d-3}\right),
		\left(\frac{4p}{p+2},\frac{2dp}{(d-5)p+2}\right),
		\left(\frac{2p}{p-2},\frac{2d}{d-2}\right)
		\right\}.
	\end{equation}
	Since $\mathcal A_Y\subset\mathcal A_d$ and
	$\frac{2d}{d-2}\in\mathcal M_0$,
	Propositions~\ref{P:STE2} and~\ref{P:STE} imply that
	\[
	\|F^\omega\|_{Y(\mathbb R)}<\infty
	\qquad\text{almost surely},
	\]
	provided that $p$ is sufficiently large.
	Fix such an $\omega$ and set
$
	F=F^\omega=e^{it\Delta}f^\omega.
$
	By the absolute continuity of the $Y$-norm, there exists an interval
	$I\ni0$ on which $\|F\|_{Y(I)}$ is sufficiently small. Proposition
	\ref{P:lwp_hartree} then yields a unique solution $v$ to
	\eqref{NLS-hartree-v1} with $v(0)=0$. Hence
$
	u=e^{it\Delta}f^\omega+v
$
	solves \eqref{HE1} with $u(0)=f^\omega$, proving the almost sure local
	well-posedness.
	
	Now suppose that
$
	\sup_{t\in I_{\max}}\|v(t)\|_{H_x^1}<\infty.
$
	Applying Proposition~\ref{P:H1-implies-scattering-hartree} forward in
	time and its time-reversed analogue backward in time, we obtain
	$I_{\max}=\mathbb R$ and $v_\pm\in H^1(\mathbb R^d)$ such that
	\[
	\|v(t)-S(t)v_\pm\|_{H_x^1}\to0
	\qquad\text{as }t\to\pm\infty.
	\]
	Therefore,
	\[
	\bigl\|u(t)-e^{it\Delta}(f^\omega+v_\pm)\bigr\|_{H_x^1}
	=
	\|v(t)-S(t)v_\pm\|_{H_x^1}
	\to0
	\]
	as $t\to\pm\infty$, which proves the claimed scattering.
\end{proof}

\section{Weighted and interaction Morawetz estimates}\label{Morawetz}

In this section we prove weighted and interaction Morawetz estimates for the
solution $v$ to \eqref{NLS-hartree-v1}. The equation \eqref{NLS-hartree-v1}
can be written in the perturbed form
\begin{equation}\label{eq:perturbed}
	(i\partial_t+\Delta)v=\mathcal H(v)+\mathcal N,
\end{equation}
where
$
\mathcal N
:=
\bigl(|x|^{-4}*|v+F|^2\bigr)(v+F)
-
\bigl(|x|^{-4}*|v|^2\bigr)v
$ and $F$ solves \eqref{eq:F-free}. 

To control the perturbative term $\mathcal N$, we
introduce the following spacetime norm. Define
$
\mathcal R_1
:=
\left\{
\frac{2d}{d-3},
\frac{2d}{d-4},
\frac{2d}{d-5}
\right\}
$
and
$
\mathcal R_2
:=
\left\{
\frac{4d}{2d-3},
\frac{2d}{d-2},
\frac{4d}{2d-5}
\right\}.
$
Then set
\begin{equation}\label{def:X}
	\begin{aligned}
		\|F\|_{\mathcal X(I)}
		:={}&
		\left\|\langle x\rangle^{\frac3p}F\right\|_
		{L_t^{\frac{2p}{p-2}}L_x^{\frac{2d}{d-5}}}
		+
		\sum_{m\in\mathcal M_0}
		\left\|\langle x\rangle^{\frac3p}\nabla F\right\|_
		{L_t^{\frac{2p}{p-2}}L_x^m}
+
		\sum_{r\in\mathcal R_1}
		\|F\|_{L_t^2L_x^r}
\\ &\quad+
		\sum_{r\in\mathcal R_2}
		\|F\|_{L_t^4L_x^r}
		+
		\|F\|_{L_t^{\frac{p}{2}}L_x^{\frac{2d}{d-3}}}+
		\|F\|_{L_t^pL_x^{\frac{2d}{d-2}}},
	\end{aligned}
\end{equation}
where $\mathcal M_0$ is defined by \eqref{def:M0} and all spacetime norms are taken over $I\times\mathbb R^d$, $d\ge5$, and
$p$ is a fixed sufficiently large finite parameter. For $d=5$, the norm involving $\frac{2d}{d-5}$ is replaced by the corresponding $L_x^\infty$ norm.

For the free evolution $F^\omega=e^{it\Delta}f^\omega$,  Propositions~\ref{P:STE2} and \ref{P:STE} imply that, if
$f\in H_{\mathrm{rad}}^s(\mathbb R^d)$ with
$s>\frac{d}{2(d-1)}$, then
$$
\|F^\omega\|_{\mathcal X(\mathbb R)}<\infty
\quad\text{almost surely},
$$
provided $p$ is chosen sufficiently large.

Throughout, we assume that the solution $v:[0,T]\times\R^d\to\C$ to \eqref{NLS-hartree-v1} satisfies the uniform mass bound
\begin{equation}\label{mass-bd}
	\|v\|_{L_t^\infty L_x^2([0,T]\times\mathbb R^d)}\lesssim 1,
\end{equation}
and define
\begin{equation}\label{bold-E} 
	\EE(T):=\sup_{t\in[0,T]}E[v(t)], 
\end{equation} 
where the energy $E[\,\cdot\,]$ is as in \eqref{E}.

We begin by recalling the standard Morawetz identity. Let $a:\mathbb{R}^d\to\mathbb{R}$ be a smooth real-valued weight centered at zero, and define
\[
M_a^0(t) := 2\Im\int_{\mathbb{R}^d} a_k(x)\,v_k(t,x)\,\bar v(t,x)\,dx
= 2\Im\int_{\mathbb{R}^d}\nabla a\cdot\nabla v\,\bar v\,dx,
\]
where $a_k=\partial_k a$, $v_k=\partial_k v$ and repeated indices are summed. The following lemma can be found in
\citep[Lemma 6.1]{Miao11}.
\begin{lem}\label{lem:3.2}
	Let $v$ solve \eqref{eq:perturbed}. Then
	\begin{align}
		\partial_t M_a^0(t)
		&=
		\int_{\mathbb R^d}
		\bigl[-\Delta\Delta a\,|v|^2
		+
		4a_{jk}\operatorname{Re}(\partial_j\overline v\,\partial_k v)
		\bigr]\,dx
	\label{eq:Ma0-L}
	\\
		&\quad
		+
		2\operatorname{Re}\int_{\mathbb R^d}\Delta a\,\overline v\,\mathcal H(v)\,dx
		+
		4\operatorname{Re}\int_{\mathbb R^d}\nabla a\cdot\nabla\overline v\,\mathcal H(v)\,dx
		\label{eq:Ma0-H}
		\\
		&\quad
		+
		2\operatorname{Re}\int_{\mathbb R^d}\Delta a\,\overline v\,\mathcal N\,dx
		+
		4\operatorname{Re}\int_{\mathbb R^d}\nabla a\cdot\nabla\overline v\,\mathcal N\,dx .
		\label{eq:Ma0-N}
	\end{align}
\end{lem}

Using Lemma~\ref{lem:3.2}, we derive the weighted Morawetz estimate for the solution $v$ of \eqref{eq:perturbed}.

\begin{prop}\label{prop:3.3}
	Let $d\ge5$ and let $v$ be the solution of \eqref{eq:perturbed} satisfying \eqref{mass-bd}. Then for any interval
	$[0,T]\subset\mathbb{R}$,
	\begin{equation}\label{E:Mor-single}
		\begin{aligned}
			&\|\langle x\rangle^{-\frac32}v\|_{L^2_{t,x}}^2
			+\|\langle x\rangle^{-\frac32}\nabla v\|_{L^2_{t,x}}^2\\
			&\quad\lesssim
			\EE(T)^{\frac{1}{2}}+
			\left|\operatorname{Re}
			\int_0^T
			\left[
			\mathfrak B_4
			\bigl(
			|v+F|^2,
			\Delta a\,(v+F)\overline v
			\bigr)(t)
			-
			\mathfrak B_4
			\bigl(
			|v|^2,
			\Delta a\,|v|^2
			\bigr)(t)
			\right]dt
			\right|\\
			&\qquad+
			\left|\operatorname{Re}
			\int_0^T
			\left[
			\mathfrak B_4
			\bigl(
			|v+F|^2,
			(v+F)\nabla a\cdot\nabla\overline v
			\bigr)(t)
			-
			\mathfrak B_4
			\bigl(
			|v|^2,
			v\nabla a\cdot\nabla\overline v
			\bigr)(t)
			\right]dt
			\right|.
		\end{aligned}
	\end{equation}
	where all spacetime norms are taken over $[0,T]\times\mathbb{R}^d$.
\end{prop}
\begin{proof}
	
	We choose the smooth radial weight $a(x):=\langle x\rangle=(1+|x|^2)^{1/2}$. A direct computation gives
\begin{equation}\label{eq:weight-a-identities}
	\begin{aligned}
		\nabla a
		&=\frac{x}{\langle x\rangle},
		\qquad
		a_{jk}
		=\frac{\delta_{jk}}{\langle x\rangle}
		-\frac{x_jx_k}{\langle x\rangle^{3}},
		\\
		\Delta a
		&=\frac{d-1}{\langle x\rangle}
		+\frac{1}{\langle x\rangle^{3}},
		\qquad
		-\Delta\Delta a
		=
		\frac{(d-1)(d-3)}{\langle x\rangle^{3}}
		+\frac{6(d-3)}{\langle x\rangle^{5}}
		+\frac{15}{\langle x\rangle^{7}} .
	\end{aligned}
\end{equation}
	
	We first bound the linear terms \eqref{eq:Ma0-L}. Noting that
$$
		\int_{\mathbb{R}^d}(-\Delta\Delta a)\,|v|^2\,dx
		\ \ge\ (d-1)(d-3)\int_{\mathbb{R}^d}\frac{|v|^2}{\langle x\rangle^{3}}\,dx,
$$
and
$$
		4\Re\int_{\mathbb{R}^d}a_{jk}\,\bar v_j v_k\,dx
		\ \ge\ 4\int_{\mathbb{R}^d}\frac{|\nabla v|^2}{\langle x\rangle^{3}}\,dx .
$$
	Then, we show that the terms in \eqref{eq:Ma0-H} are nonnegative.   Since
$
	\nabla\bigl(|x|^{-4}\ast |v|^2\bigr)
	=
	(\nabla |x|^{-4})\ast |v|^2
$
	and $\nabla |x|^{-4}$ is odd, an integration by parts followed by
	symmetrisation in $x$ and $y$ gives
\begin{equation*}
	\begin{aligned}
		&2\Re\!\int_{\mathbb R^d}\!\Delta a\,\bar v\,
		\bigl(|x|^{-4}\ast |v|^2\bigr)v\,dx
		+4\Re\!\int_{\mathbb R^d}\!\nabla a\cdot\nabla\bar v\,
		\bigl(|x|^{-4}\ast |v|^2\bigr)v\,dx  \\
		&\qquad
		=-2\int_{\mathbb R^d}
		\nabla a\cdot\nabla\bigl(|x|^{-4}\ast |v|^2\bigr)\,
		|v|^2\,dx  \\
		&\qquad
		=4\iint_{\mathbb R^d\times\mathbb R^d}
		\frac{\bigl(\nabla a(x)-\nabla a(y)\bigr)\cdot(x-y)}
		{|x-y|^6}\,
		|v(x)|^2|v(y)|^2\,dx\,dy\ge 0 .
	\end{aligned}
\end{equation*}	where the inequality follows from the convexity of $a(x)=\langle x\rangle$.
		
 By Cauchy--Schwarz,
$|\nabla a|\le1$, and the mass bound \eqref{mass-bd}, 
\[ \sup_{t\in[0,T]}|M_a^0(t)| \leq 2\|\nabla a\|_{L^\infty} \|v(t)\|_{L_t^{\infty}L_x^2}\|\nabla v(t)\|_{L_t^{\infty}L_x^2} \lesssim \EE(T)^{1/2}. \]
Then, integrating the identity in Lemma~\ref{lem:3.2} over $[0,T]$ and using \eqref{eq:weight-a-identities}, we obtain \eqref{E:Mor-single}.
\end{proof}

We next turn to the interaction Morawetz estimate.  Let
$
m(x-y)=|x-y|.
$
We define the interaction Morawetz action by
\begin{equation}\label{def:interaction-M}
	M(t):=
	2\Im\iint_{\mathbb R^d\times\mathbb R^d}
	|v(y)|^{2}\,\nabla m(x-y)\cdot\nabla v(x)\,\overline v(x)\,dxdy.
\end{equation}
For later use, we also introduce the momentum bracket
\[
\{f,g\}_p^j
:=
\Re\bigl(f\partial_j\overline g-g\partial_j\overline f\bigr),
\qquad f,g:\mathbb R^d\to\mathbb C.
\]
Differentiating $M(t)$, using the mass identity
\[
\partial_t |v|^2
=
-2\partial_j\Im(\partial_j v\,\overline v)
+
2\Im(\mathcal N\overline v),
\]
where $\Im(\mathcal H(v)\overline v)\equiv0$, and applying the preceding
Morawetz identity with spatial weight $m(x-y)$, we obtain the
following interaction identity.
 \begin{lem}\label{lem:interaction-id}
 	Let $v$ solve \eqref{eq:perturbed}. Then, for every smooth real-valued weight $m$,
 	\begin{align}
 		\dot M(t)
 		&=\iint(-\Delta\Delta m)(x-y)\,|v(x)|^{2}|v(y)|^{2}\,dxdy
 		\label{eq:int-lin}\\
 		&\quad+4\iint|v(y)|^{2}m_{jk}(x-y)\,\Re(\partial_j\bar v\,\partial_k v)(x)\,dxdy
 		\label{eq:int-curr1}\\
 		&\quad-4\iint\Im(\partial_jv\,\bar v)(y)\,m_{jk}(x-y)\,\Im(\partial_kv\,\bar v)(x)\,dxdy
 		\label{eq:int-curr2}\\
 		&\quad+4\iint\Im(\mathcal N\bar v)(y)\,\nabla m(x-y)\cdot\Im(\nabla v\,\bar v)(x)\,dxdy
 		\label{eq:int-Nsrc}\\
 		&\quad+2\iint|v(y)|^{2}m_j(x-y)\,\{\mathcal H(v)+\mathcal N,v\}_p^{j}(x)\,dxdy.
 		\label{eq:int-brk}
 	\end{align}
 \end{lem}
 
 Using Lemma~\ref{lem:interaction-id}, we derive the interaction Morawetz estimate.

 \begin{prop}\label{prop:interaction}
 	Let $d\ge5$ and let $v$ solve \eqref{eq:perturbed} and satisfy
 	\eqref{mass-bd}. Then, for every interval $[0,T]\subset\mathbb R$,
 	\begin{equation}\label{E:Mor-interaction}
 		\begin{aligned}
 			\bigl\||\nabla|^{-\frac{d-3}{2}}|v|^2
 			\bigr\|_{L^2_{t,x}}^2
 			\lesssim{}&
 			\EE(T)^{\frac12}+\EE(T)^{\frac12}
 			\int_0^T
 			\mathfrak B_4
 			\bigl(|v+F|^2,|\Im(F\overline v)|\bigr)(t)\,dt\\
 			&+\left|
 			\int_0^T
 			\Bigl[
 			\mathfrak B_4
 			\Bigl(
 			2\Re(v\overline F)+|F|^2,
 			\bigl(|x|^{-1}*|v|^2\bigr)|v|^2
 			\Bigr)(t)\right.\\
 			&\left.
 			+\mathfrak B_4
 			\Bigl(
 			|v+F|^2,
 			\bigl(|x|^{-1}*|v|^2\bigr)\Re(F\overline v)
 			\Bigr)(t)
 			\Bigr]\,dt
 			\right|\\
 			&+\left|
 			\int_0^T\!\!\iint_{\mathbb R^{2d}}
 			|v(t,y)|^2m_j(x-y)
 			\Re\bigl(
 			\mathcal N(t,x)\partial_j\overline v(t,x)
 			\bigr)\,dx\,dy\,dt
 			\right|.
 		\end{aligned}
 	\end{equation}
 \end{prop}

\begin{proof}
	We take $m(x-y)=|x-y|$. Then, for $d\ge 5$, 
	\begin{equation}\label{eq:weight-m-id}
		\begin{gathered}
			\nabla m(x-y)=\frac{x-y}{|x-y|},\qquad
			m_{jk}(x-y)=\frac{\delta_{jk}}{|x-y|}-\frac{(x_j-y_j)(x_k-y_k)}{|x-y|^3},\\[2pt]
			\Delta m(x-y)=\frac{d-1}{|x-y|},\qquad
			-\Delta\Delta m(x-y)=\frac{(d-1)(d-3)}{|x-y|^3}.
		\end{gathered}
	\end{equation}
	Substituting \eqref{eq:weight-m-id} into \eqref{eq:int-lin}, we obtain
	\begin{equation}\label{eq:int-lin-eval}
		\eqref{eq:int-lin}=(d-1)(d-3)\iint\frac{|v(x)|^2\,|v(y)|^2}{|x-y|^{3}}\,dxdy.
	\end{equation}
	
	We claim that
	\begin{equation}\label{eq:claim-curr}
	\eqref{eq:int-curr1}	+\eqref{eq:int-curr2}\ \ge\ 0.
	\end{equation}
	Indeed, since the matrix $(m_{jk})$ is positive semidefinite, the pointwise inequality
	\[
	|v(y)|^2\,m_{jk}(x-y)\,\Re(\partial_j\bar v\,\partial_k v)(x)
	\ \ge\ 
	|v(y)|^2\,m_{jk}(x-y)\,\frac{\Im(\bar v\,\partial_j v)\,\Im(\bar v\,\partial_k v)}{|v|^2}(x)
	\]
	holds. Symmetrizing in \eqref{eq:int-curr1} in the variables $x$ and $y$ and applying the Cauchy--Schwarz inequality associated with the quadratic form $(m_{jk})$, we obtain \eqref{eq:claim-curr}.
	
	Next, we decompose the bracket term in \eqref{eq:int-brk} as
$
	\{\mathcal H(v)+\mathcal N,v\}_p^{j}=\{\mathcal H(v),v\}_p^{j}+\{\mathcal N,v\}_p^{j}.
$
 Using
	$\{\mathcal H(v),v\}_p^{j}=-|v|^2\,\partial_j\bigl(|x|^{-4}\ast|v|^2\bigr)$ and
	$\partial_j\bigl(|x|^{-4}\ast|v|^2\bigr)(x)=-4\!\int\!\frac{(x-z)_j}{|x-z|^{6}}|v(z)|^2\,dz$,
symmetrisation in $x$ and $z$ yields
	\begin{equation}\label{eq:Hartree-pos}
		\begin{aligned}
			&2\iint|v(y)|^{2}\,m_j(x-y)\,\{\mathcal H(v),v\}_p^{j}(x)\,dxdy\\
			=&8\iiint\frac{m_j(x-y)\,(x-z)_j}{|x-z|^{6}}\,|v(x)|^2|v(y)|^2|v(z)|^2\,dxdydz\\
			=&4\iiint\frac{\bigl[\nabla m(x-y)-\nabla m(z-y)\bigr]\cdot(x-z)}{|x-z|^{6}}\,|v(x)|^2|v(y)|^2|v(z)|^2\,dxdydz\ \ge\ 0,
		\end{aligned}
	\end{equation}
	where the last inequality follows from the convexity of \(m=|\cdot|\).
	
	Accordingly, using the results of 
	\eqref{eq:claim-curr} and \eqref{eq:Hartree-pos}, integrating the
	identity in Lemma~\ref{lem:interaction-id} over $[0,T]$ yields
	\begin{equation}\label{eq:after-FTC}
		(d-1)(d-3)
		\iint_{[0,T]\times\R^{2d}}
		\frac{|v(t,x)|^2|v(t,y)|^2}{|x-y|^3}
		\,dx\,dy\,dt
		\le
		2\sup_{t\in[0,T]}|M(t)|
		+\mathcal I_1+\mathcal I_2,
	\end{equation}
	where
	\[
	\mathcal I_1
	:=
	\left|\int_0^T \eqref{eq:int-Nsrc}\,dt\right|,
	\qquad
	\mathcal I_2
	:=
	\left|
	\int_0^T\!\!\iint
	2|v(t,y)|^2m_j(x-y)
	\{\mathcal N,v\}_p^j(t,x)
	\,dx\,dy\,dt
	\right|.
	\]
	
	We next estimate $\mathcal I_1$ and $\mathcal I_2$. For $\mathcal I_1$,
	using $|\nabla m|\le 1$, H\"older's inequality, and \eqref{mass-bd},
	we obtain
	\begin{align}
		\mathcal I_1
		&\le
		4\int_0^T\!\!\iint_{\mathbb R^d\times\mathbb R^d}
		\bigl|\Im(\mathcal N(t,y)\overline{v(t,y)})\bigr|
		|\nabla v(t,x)|\,|v(t,x)|
		\,dx\,dy\,dt
		\notag\\
		&\le
		4\int_0^T
		\left(
		\int_{\mathbb R^d}
		\bigl|\Im(\mathcal N(t,y)\overline{v(t,y)})\bigr|
		\,dy
		\right)
		\|v(t)\|_{L_x^2}
		\|\nabla v(t)\|_{L_x^2}
		\,dt
		\notag\\
		&\lesssim
			\EE(T)^{\frac12}
		\int_0^T
		\mathfrak B_4
		\bigl(|v+F|^2,|\Im(F\overline v)|\bigr)(t)\,dt.
		\label{eq:I1-Hartree-form}
	\end{align}

	Integrating by parts in $x_j$ and using
	\eqref{eq:weight-m-id}, we have
\begin{align}
	\mathcal I_2
	\lesssim{}&
	\left|
	\int_0^T\!\!\iint_{\mathbb R^{2d}}
	|v(t,y)|^2m_j(x-y)
	\Re\bigl(
	\mathcal N(t,x)\partial_j\overline v(t,x)
	\bigr)
	\,dx\,dy\,dt
	\right|
	\notag\\
	&+
	\left|
	\int_0^T\int_{\mathbb R^d}
	\bigl(|x|^{-1}*|v(t)|^2\bigr)(x)
	\Re\bigl(
	\mathcal N(t,x)\overline v(t,x)
	\bigr)
	\,dx\,dt
	\right|	\notag\\
	\lesssim{}&
	\left|
	\int_0^T\!\!\iint_{\mathbb R^{2d}}
	|v(t,y)|^2m_j(x-y)
	\Re\bigl(
	\mathcal N(t,x)\partial_j\overline v(t,x)
	\bigr)\,dx\,dy\,dt
	\right|
	\notag\\
	&+
	\left|
	\int_0^T
	\Bigl[
	\mathfrak B_4
	\Bigl(
	2\Re(v\overline F)+|F|^2,
	\bigl(|x|^{-1}*|v|^2\bigr)|v|^2
	\Bigr)(t)
	+
	\mathfrak B_4
	\Bigl(
	|v+F|^2,
	\bigl(|x|^{-1}*|v|^2\bigr)\Re(F\overline v)
	\Bigr)(t)
	\Bigr]\,dt
	\right|.
	\label{eq:I2-final}
\end{align}

Finally,
\begin{equation}\label{eq:M-bound-interaction}
	\sup_{t\in[0,T]}|M(t)|
	\lesssim
	\sup_{t\in[0,T]}
	\|v(t)\|_{L_x^2}^3\|\nabla v(t)\|_{L_x^2}
	\lesssim
	\EE(T)^{\frac12}.
\end{equation}
Substituting \eqref{eq:I1-Hartree-form}, \eqref{eq:I2-final}, and
\eqref{eq:M-bound-interaction} into \eqref{eq:after-FTC}, and using Lemma~\ref{lem:2.3} with
$\alpha=(d-3)/2$ gives \eqref{E:Mor-interaction}. This completes the
proof.
\end{proof}

By combining the weighted and interaction Morawetz estimates with Lemma~\ref{lem:Hartree-reduction} in the case \(\gamma=4\), we obtain the following spacetime bound.

\begin{prop}\label{P:Morawetz}
	Let $d\ge5$, and let $v$ be the solution of \eqref{eq:perturbed}
	on $[0,T]$ satisfying \eqref{mass-bd}, with $F$ denoting the associated
	free evolution. Define
	\begin{equation}\label{def:A}
		\begin{aligned}
			A(T)
			:=&
			\|\x^{-\frac32}v\|_{L^2_{t,x}([0,T]\times\R^d)}^2
			+
			\|\x^{-\frac32}\nabla v\|_{L^2_{t,x}([0,T]\times\R^d)}^2
			+
			\left\||\nabla|^{-\frac{d-3}{2}}|v|^2\right\|_{L^2_{t,x}([0,T]\times\R^d)}^2 .
		\end{aligned}
	\end{equation}
	Then
	\begin{equation}\label{E:Morawetz2}
		\begin{aligned}
		A(T)
		\lesssim\;&
		\EE(T)^{\frac12}
		+\EE(T)
		\|F\|_{\mathcal X([0,T])}^{\frac{2p}{p-2}}+
		\EE(T)^{\frac{p-2}{2(p-1)}}
		\left(
		\|F\|_{\mathcal X([0,T])}^{\frac{2p}{p-1}}
		+
		\|F\|_{\mathcal X([0,T])}^{\frac{3p}{p-1}}
		\right)\\
		&+
		\|F\|_{\mathcal X([0,T])}^2
		+\|F\|_{\mathcal X([0,T])}^4+
		\EE(T)
		\left(
		\|F\|_{\mathcal X([0,T])}^2
		+\|F\|_{\mathcal X([0,T])}^3
		+\|F\|_{\mathcal X([0,T])}^4
		\right).
		\end{aligned}
	\end{equation}
	Here $\EE(T)$ is defined in \eqref{bold-E}.
\end{prop}
\begin{proof}
Combining the definition of $A(T)$ with the weighted Morawetz estimate
\eqref{E:Mor-single} and the interaction Morawetz estimate
\eqref{E:Mor-interaction}, we obtain 
\begin{subequations}\label{E:Mor}
	\begin{align}
		A(T)
		\lesssim{}&
		\EE(T)^{\frac{1}{2}}
		\label{E:Mor-a}\\
		&+
		\left|\operatorname{Re}
		\int_0^T
		\left[
		\mathfrak B_4
		\bigl(
		|v+F|^2,
		\Delta a\,(v+F)\overline v
		\bigr)(t)
		-
		\mathfrak B_4
		\bigl(
		|v|^2,
		\Delta a\,|v|^2
		\bigr)(t)
		\right]dt
		\right|
		\label{E:Mor-b}\\
		&+
		\left|\operatorname{Re}
		\int_0^T
		\left[
		\mathfrak B_4
		\bigl(
		|v+F|^2,
		(v+F)\nabla a\cdot\nabla\overline v
		\bigr)(t)
		-
		\mathfrak B_4
		\bigl(
		|v|^2,
		v\nabla a\cdot\nabla\overline v
		\bigr)(t)
		\right]dt
		\right|
		\label{E:Mor-c}\\
		&+
		\EE(T)^{\frac12}
		\int_0^T
		\mathfrak B_4
		\bigl(
		|v+F|^2,
		|\Im(F\overline v)|
		\bigr)(t)\,dt
		\label{E:Mor-d}\\
		&+
		\left|
		\int_0^T
		\left[
		\mathfrak B_4
		\Bigl(
		2\Re(v\overline F)+|F|^2,
		\bigl(|x|^{-1}*|v|^2\bigr)|v|^2
		\Bigr)(t)
		+
		\mathfrak B_4
		\Bigl(
		|v+F|^2,
		\bigl(|x|^{-1}*|v|^2\bigr)\Re(F\overline v)
		\Bigr)(t)
		\right]dt
		\right|
		\label{E:Mor-e}\\
		&+
		\left|
		\int_0^T\!\!\iint_{\mathbb R^{2d}}
		|v(t,y)|^2m_j(x-y)
		\Re\bigl(
		\mathcal N(t,x)\partial_j\overline v(t,x)
		\bigr)\,dx\,dy\,dt
		\right|.
		\label{E:Mor-f}
	\end{align}
\end{subequations}

	We shall use the following interpolation estimates. By H\"older's inequality, Sobolev embedding, and the definition of $A(T)$, for $2\le p\le\infty$, we have	\begin{align}
		\|\x^{-\frac3p}v\|_{L_t^p L_x^{\frac{2pd}{pd-2p+4}}}
		&\lesssim\|v\|_{L_t^\infty L_x^{\frac{2d}{d-2}}}^{1-\frac2p}
		\|\x^{-\frac32}v\|_{L_{t,x}^2}^{\frac2p}
		\lesssim\EE(T)^{\frac12-\frac1p}A(T)^{\frac1p},\label{interpolate1}\\
		\|\x^{-\frac3p}\nabla v\|_{L_t^p L_x^2}
		&\lesssim\|\nabla v\|_{L_t^\infty L_x^2}^{1-\frac2p}
		\|\x^{-\frac32}\nabla v\|_{L_{t,x}^2}^{\frac2p}
		\lesssim\EE(T)^{\frac12-\frac1p}A(T)^{\frac1p},\label{interpolate2}
	\end{align}
	
Taking $\gamma=4$, $a=1$, and $b=b'=2$ in
Lemma~\ref{lem:Hartree-reduction}, with
$s=\frac{d-3}{2}$ and $s=\frac{d-5}{2}$ for the first and second terms,
respectively, and applying Lemma~\ref{lem:1}, we obtain
\begin{equation}\label{eq:Hartree-pairing-consequence}
	\|\mathfrak B_4(|v|^2,h)\|_{L_t^1}
	+\|\mathfrak B_4(h,|v|^2)\|_{L_t^1}
	\lesssim
	A(T)^{\frac12}
	\|h\|_{L_t^2L_x^{\frac{2d}{2d-5}}}.
\end{equation}
Then, combining \eqref{mass-bd} and \eqref{eq:Hartree-pairing-consequence}
with H\"older's inequality, Sobolev embedding, and Lemma~\ref{lem:1},
we have
\begin{align}
	\eqref{E:Mor-b}
	\lesssim{}&
	\left\|
	\mathfrak B_4
	\bigl(
	2\Re(v\overline F)+|F|^2,
	\Delta a\,|v|^2
	\bigr)
	\right\|_{L_t^1}+
	\left\|
	\mathfrak B_4
	\bigl(
	|v|^2,
	\Delta a\,F\overline v
	\bigr)
	\right\|_{L_t^1}+
	\left\|
	\mathfrak B_4
	\bigl(
	2\Re(v\overline F)+|F|^2,
	\Delta a\,F\overline v
	\bigr)
	\right\|_{L_t^1}
	\notag\\
	\lesssim{}&
	A(T)^{\frac12}
	\left(
	\|vF\|_{L_t^2L_x^{\frac{2d}{2d-5}}}
	+
	\|F^2\|_{L_t^2L_x^{\frac{2d}{2d-5}}}
	\right)
	+
	\left(
	\|vF\|_{L_t^2L_x^{\frac{2d}{2d-3}}}
	+
	\|F^2\|_{L_t^2L_x^{\frac{2d}{2d-3}}}
	\right)
	\|Fv\|_{L_t^2L_x^{\frac{2d}{2d-5}}}
	\notag\\
	\lesssim{}&
	A(T)^{\frac12}
	\left(
	\|F\|_{L_t^2L_x^{\frac{2d}{d-5}}}
	\|v\|_{L_t^\infty L_x^2}
	+
	\|F\|_{L_t^4L_x^{\frac{4d}{2d-5}}}^2
	\right)
	+
	\left(
	\|F\|_{L_t^2L_x^{\frac{2d}{d-3}}}
	\|v\|_{L_t^\infty L_x^2}
	+
	\|F\|_{L_t^4L_x^{\frac{4d}{2d-3}}}^2
	\right)
	\|F\|_{L_t^2L_x^{\frac{2d}{d-3}}}
	\|v\|_{L_t^\infty L_x^{\frac{2d}{d-2}}}
	\label{eq:Mor-b-intermediate}	\notag\\
	\lesssim{}&
	A(T)^{\frac12}
	\left(
	\|F\|_{\mathcal X([0,T])}
	+
	\|F\|_{\mathcal X([0,T])}^2
	\right)
	+
	\EE(T)^{\frac12}
	\left(
	\|F\|_{\mathcal X([0,T])}^2
	+
	\|F\|_{\mathcal X([0,T])}^3
	\right).
\end{align}	
	
To estimate the term in \eqref{E:Mor-c}, we integrate by parts to obtain
	\begin{align}
		&\Re\,
		\mathfrak B_4
		\bigl(
		2\Re(v\overline F)+|F|^2,
		v\nabla a\cdot\nabla\overline v
		\bigr)\notag\\
		&=
		-\frac12
		\mathfrak B_4
		\bigl(
		\nabla(2\Re(v\overline F)+|F|^2),
		|v|^2\nabla a
		\bigr)
		-\frac12
		\mathfrak B_4
		\bigl(
		2\Re(v\overline F)+|F|^2,
		|v|^2\Delta a
		\bigr).
		\label{eq:weighted-integration-by-parts}
	\end{align}
	Using \eqref{eq:weight-a-identities} and
	\eqref{eq:Hartree-pairing-consequence}, we have
\begin{align}
	\eqref{E:Mor-c}
	\lesssim{}&
	\left\|
	\mathfrak B_4
	\bigl(
	\nabla(2\Re(v\overline F)+|F|^2),
	|v|^2\nabla a
	\bigr)
	\right\|_{L_t^1}
	+
	\left\|
	\mathfrak B_4
	\bigl(
	2\Re(v\overline F)+|F|^2,
	|v|^2\Delta a
	\bigr)
	\right\|_{L_t^1}
	\notag\\
	&\quad+
	\left\|
	\mathfrak B_4
	\bigl(
	|v|^2,
	F\nabla a\cdot\nabla\overline v
	\bigr)
	\right\|_{L_t^1}
	+
	\left\|
	\mathfrak B_4
	\bigl(
	2\Re(v\overline F)+|F|^2,
	F\nabla a\cdot\nabla\overline v
	\bigr)
	\right\|_{L_t^1}.
	\label{eq:weight-a-identities1}
\end{align}
Applying H\"older's inequality and
		\eqref{interpolate1}--\eqref{interpolate2}, we obtain
			\begin{align} 
		\mathrm{RHS}\eqref{eq:weight-a-identities1}	&\lesssim
		A(T)^{\frac12}
		\Bigl(
		\|F\nabla v\|_{L_t^2L_x^{\frac{2d}{2d-5}}}
		+\|v\nabla F\|_{L_t^2L_x^{\frac{2d}{2d-5}}}
		+\|F\nabla F\|_{L_t^2L_x^{\frac{2d}{2d-5}}}
		+\|vF\|_{L_t^2L_x^{\frac{2d}{2d-5}}}
		+\|F^2\|_{L_t^2L_x^{\frac{2d}{2d-5}}}
		\Bigr)
		\notag\\
		&\quad+
		\left(
		\|F\|_{L_t^2L_x^{\frac{2d}{d-3}}}
		\|v\|_{L_t^\infty L_x^2}
		+
		\|F\|_{L_t^4L_x^{\frac{4d}{2d-3}}}^2
		\right)
		\|\x^{\frac3p}F\|_
		{L_t^{\frac{2p}{p-2}}L_x^{\frac{2d}{d-5}}}
		\|\x^{-\frac3p}\nabla v\|_{L_t^pL_x^2}
		\notag\\
		&\lesssim
		A(T)^{\frac12}
		\Biggl(
		\|\x^{\frac3p}F\|_
		{L_t^{\frac{2p}{p-2}}L_x^{\frac{2d}{d-5}}}
		\|\x^{-\frac3p}\nabla v\|_{L_t^pL_x^2}
		+
		\|\x^{-\frac3p}v\|_
		{L_t^pL_x^{\frac{2pd}{pd-2p+4}}}
		\|\x^{\frac3p}\nabla F\|_
		{L_t^{\frac{2p}{p-2}}
			L_x^{\frac{2pd}{pd-3p-4}}}
		\notag\\
		&\quad+
		\|F\|_{L_t^pL_x^{\frac{2d}{d-2}}}
		\|\x^{\frac3p}\nabla F\|_
		{L_t^{\frac{2p}{p-2}}L_x^{\frac{2d}{d-3}}}
		+
		\|v\|_{L_t^\infty L_x^2}
		\|F\|_{L_t^2L_x^{\frac{2d}{d-5}}}
		+
		\|F\|_{L_t^4L_x^{\frac{4d}{2d-5}}}^2
		\Biggr)
		\notag\\
			&\quad+
		\left(
		\|F\|_{L_t^2L_x^{\frac{2d}{d-3}}}
		\|v\|_{L_t^\infty L_x^2}
		+
		\|F\|_{L_t^4L_x^{\frac{4d}{2d-3}}}^2
		\right)
		\|\x^{\frac3p}F\|_
		{L_t^{\frac{2p}{p-2}}L_x^{\frac{2d}{d-5}}}
		\|\x^{-\frac3p}\nabla v\|_{L_t^pL_x^2}
		\notag\\
		&\lesssim
		\EE(T)^{\frac12-\frac1p}
		A(T)^{\frac12+\frac1p}
		\|F\|_{\mathcal X([0,T])}
		+
		A(T)^{\frac12}
		\left(
		\|F\|_{\mathcal X([0,T])}
		+\|F\|_{\mathcal X([0,T])}^2
		\right)
		\notag\\
		&\qquad+
		\EE(T)^{\frac12-\frac1p}
		A(T)^{\frac1p}
		\left(
		\|F\|_{\mathcal X([0,T])}^2
		+\|F\|_{\mathcal X([0,T])}^3
		\right).
		\label{eq:combined-weighted-pairing-bounds}
	\end{align}
The following terms are estimated similarly. Using the same choice of $s$ and $b$, H\"older's inequality, and \eqref{interpolate1}--\eqref{interpolate2}, we have:
\begin{align}
	\eqref{E:Mor-d}
	\lesssim{}&
	\EE(T)^{\frac12}A(T)^{\frac12}
	\|F\|_{\mathcal X([0,T])}
	+
	\EE(T)
	\left(
	\|F\|_{\mathcal X([0,T])}^2
	+
	\|F\|_{\mathcal X([0,T])}^3
	\right).
	\label{eq:Mor-d-bound}
\end{align}

	We now estimate the terms in \eqref{E:Mor-e} and
\eqref{E:Mor-f}. By Lemma \ref{lem:Hardy} and
\eqref{mass-bd},
\begin{align}
	\eqref{E:Mor-e}
	\lesssim{}&
	\bigl\||x|^{-1}*|v|^2\bigr\|_{L_t^\infty L_x^\infty}
	\Bigl[
	\bigl\|
	\mathfrak B_4
	\bigl(
	|vF|+|F|^2,
	|v|^2
	\bigr)
	\bigr\|_{L_t^1}
	+
	\bigl\|
	\mathfrak B_4
	\bigl(
	|v+F|^2,
	|Fv|
	\bigr)
	\bigr\|_{L_t^1}
	\notag\\
	\lesssim{}&\|v\|_{L_t^\infty L_x^2}
	\|\nabla v\|_{L_t^\infty L_x^2}
	\biggl[A(T)^{\frac12}
	\left(
	\|F\|_{\mathcal X([0,T])}
	+
	\|F\|_{\mathcal X([0,T])}^2
	\right)
+
	\EE(T)^{\frac12}
	\left(
	\|F\|_{\mathcal X([0,T])}^2
	+
	\|F\|_{\mathcal X([0,T])}^3
	\right)
	\biggr]
	\notag\\
	\lesssim{}&
	\EE(T)^{\frac12}A(T)^{\frac12}
	\left(
	\|F\|_{\mathcal X([0,T])}
	+
	\|F\|_{\mathcal X([0,T])}^2
	\right)+
	\EE(T)
	\left(
	\|F\|_{\mathcal X([0,T])}^2
	+
	\|F\|_{\mathcal X([0,T])}^3
	\right).
	\label{eq:Mor-e-bound}
\end{align}

To estimate \eqref{E:Mor-f}, we first expand
\begin{align}
	\mathcal N
	={}&
	\bigl(|x|^{-4}*|v|^2\bigr)F
	+
	\bigl(
	|x|^{-4}*
	(2\Re(v\overline F)+|F|^2)
	\bigr)v
	\notag\\
	&+
	\bigl(
	|x|^{-4}*
	(2\Re(v\overline F)+|F|^2)
	\bigr)F.
	\label{eq:N-expansion-interaction}
\end{align}
For the middle term in \eqref{eq:N-expansion-interaction}, using
$
\Re\bigl(v\partial_j\overline v\bigr)
=\frac12\partial_j|v|^2
$
and integrating by parts in $x_j$, we have
\begin{align}
	&\int_0^T\!\!\iint_{\mathbb R^{2d}}
	|v(t,y)|^2m_j(x-y)
	\bigl(
	|x|^{-4}*
	(2\Re(v\overline F)+|F|^2)
	\bigr)(t,x)
	\Re\bigl(
	v(t,x)\partial_j\overline v(t,x)
	\bigr)
	\,dx\,dy\,dt
	\notag\\
	&\quad=
	-\frac12
	\int_0^T\!\!\iint_{\mathbb R^{2d}}
	|v(t,y)|^2\Delta m(x-y)
	\bigl(
	|x|^{-4}*
	(2\Re(v\overline F)+|F|^2)
	\bigr)(t,x)
	|v(t,x)|^2
	\,dx\,dy\,dt
	\notag\\
	&\qquad
	-\frac12
	\int_0^T\!\!\iint_{\mathbb R^{2d}}
	|v(t,y)|^2m_j(x-y)
	\bigl(
	|x|^{-4}*
	\partial_j(2\Re(v\overline F)+|F|^2)
	\bigr)(t,x)
	|v(t,x)|^2
	\,dx\,dy\,dt.
	\label{eq:Mor-f-IBP}
\end{align}
By \eqref{eq:weight-m-id}, the first term on the right-hand side of
\eqref{eq:Mor-f-IBP} is estimated as follows:
\begin{align}
	&\left|
	\int_0^T
	\mathfrak B_4
	\Bigl(
	2\Re(v\overline F)+|F|^2,
	\bigl(|x|^{-1}*|v|^2\bigr)|v|^2
	\Bigr)(t)\,dt
	\right|
	\notag\\
	&\quad\le
	\int_0^T
	\int_{\mathbb R^d}
	\left|
	|x|^{-4}*
	\bigl(2\Re(v\overline F)+|F|^2\bigr)
	\right|(t,x)
	\bigl(|x|^{-1}*|v(t)|^2\bigr)(x)
	|v(t,x)|^2\,dx\,dt
	\notag\\
	&\quad\lesssim
	\bigl\||x|^{-1}*|v|^2\bigr\|_{L_t^\infty L_x^\infty}
	\left(
	\bigl\|
	\mathfrak B_4
	\bigl(
	|vF|,
	|v|^2
	\bigr)
	\bigr\|_{L_t^1}
	+
	\bigl\|
	\mathfrak B_4
	\bigl(
	|F|^2,
	|v|^2
	\bigr)
	\bigr\|_{L_t^1}
	\right)
	\notag\\
	&\quad\lesssim
	\EE(T)^{\frac12}A(T)^{\frac12}
\left(
\|F\|_{\mathcal X([0,T])}
+
\|F\|_{\mathcal X([0,T])}^2
\right).
	\label{eq:Mor-f-lower-order-first}
\end{align}

For the second term on the right-hand side of
\eqref{eq:Mor-f-IBP}, \eqref{mass-bd} and
\eqref{eq:Hartree-pairing-consequence} give
\begin{align}
	&\left|
	\int_0^T\!\!\iint_{\mathbb R^{2d}}
	|v(t,y)|^2m_j(x-y)
	\bigl(
	|x|^{-4}*
	\partial_j(2\Re(v\overline F)+|F|^2)
	\bigr)(t,x)
	|v(t,x)|^2
	\,dx\,dy\,dt
	\right|
	\notag\\
	&\qquad\lesssim
	\left\|
	\mathfrak B_4
	\bigl(
	|\nabla(2\Re(v\overline F)+|F|^2)|,
	|v|^2
	\bigr)
	\right\|_{L_t^1}
	\notag\\
	&\qquad\lesssim
		\EE(T)^{\frac12-\frac1p}
	A(T)^{\frac12+\frac1p}
	\|F\|_{\mathcal X([0,T])}
	+
	A(T)^{\frac12}
\|F\|_{\mathcal X([0,T])}^2.
	\label{eq:Mor-f-gradient}
\end{align}
The first and third terms in \eqref{eq:N-expansion-interaction}
are estimated  by
\begin{align}
	&\left\|
	\mathfrak B_4
	\bigl(
	|v|^2,
	|F\nabla v|
	\bigr)
	\right\|_{L_t^1}+\left\|
	\mathfrak B_4
	\bigl(
	|2\Re(v\overline F)+|F|^2|,
	|F\nabla v|
	\bigr)
	\right\|_{L_t^1}
	\notag\\
	&\qquad\lesssim
	\EE(T)^{\frac12-\frac1p}
	A(T)^{\frac1p}
	\left(
	\|F\|_{\mathcal X([0,T])}^2
	+\|F\|_{\mathcal X([0,T])}^3
	\right).
	\label{eq:Mor-f-higher}
\end{align}

	Combining Propositions~\ref{prop:3.3} and
	\ref{prop:interaction} with the preceding estimates yields
	\begin{align}
		A(T)
		\lesssim{}&
		\EE(T)^{\frac12}
		+\EE(T)^{\frac12-\frac1p}
		A(T)^{\frac12+\frac1p}
		\|F\|_{\mathcal X([0,T])}
+
		\EE(T)^{\frac12-\frac1p}
		A(T)^{\frac1p}
		\left(
		\|F\|_{\mathcal X([0,T])}^2
		+\|F\|_{\mathcal X([0,T])}^3
		\right)
		\notag\\
		&+
		A(T)^{\frac12}
		\left(
		\|F\|_{\mathcal X([0,T])}
		+\|F\|_{\mathcal X([0,T])}^2
		\right)
+
		\EE(T)^{\frac12}A(T)^{\frac12}
		\left(
		\|F\|_{\mathcal X([0,T])}
		+\|F\|_{\mathcal X([0,T])}^2
		\right)
		\notag\\
		&+
		\EE(T)^{\frac12}
		\left(
		\|F\|_{\mathcal X([0,T])}^2
		+\|F\|_{\mathcal X([0,T])}^3
		\right)
+
		\EE(T)
		\left(
		\|F\|_{\mathcal X([0,T])}^2
		+\|F\|_{\mathcal X([0,T])}^3
		\right).
		\label{eq:Morawetz-before-Young}
	\end{align}
		Inserting these bounds into \eqref{E:Mor} and applying Young's inequality to
	absorb the positive powers of $A(T)$ into the left-hand side yields
	\eqref{E:Morawetz2}.
\end{proof}

\section{Energy bounds for the forced equation}\label{energy-bds}

\begin{prop}[Energy bounds]\label{P:energy-bds}
Let $d\ge 7$.	Suppose that $F$ is a solution to $(i\partial_t+\Delta)F=0$ satisfying
	\begin{equation}\label{3:30}
		\|F\|_{\mathcal{X}(\R)}+\|F\|_{L_t^\infty L_x^2(\R\times\R^d)}+\|F\|_{L_t^\infty L_x^{\frac{2d}{d-2}}(\R\times\R^d)}<\infty.
	\end{equation}
	Let $v_0\in H^1(\R^d)$ and let $v:\Imax\times\R^d\to\C$ be the maximal-lifespan solution to \eqref{NLS-hartree-v1} with $v(0)=v_0$. Then $\sup_{t\in\Imax}\|v(t)\|_{H_x^1}<\infty$.
\end{prop}
\begin{proof}
	  A direct computation based on \eqref{NLS-hartree-v1} yields
	\begin{align*}
		\partial_t E[v(t)]
		&=-\Re\!\int[\partial_t\bar v]\bigl[\H(v+F)-\H(v)\bigr]dx\\
		&=-\tfrac14\partial_t\!\int\!\bigl[(|x|^{-4}\!\ast\!|v+F|^2)|v+F|^2-(|x|^{-4}\!\ast\!|v|^2)|v|^2-(|x|^{-4}\!\ast\!|F|^2)|F|^2\bigr]dx\\
		&\quad+\Re\!\int\!\bigl[\bar\H(v+F)-\bar\H(F)\bigr]\partial_t F\,dx,
	\end{align*}
Since $F$ solves $(i\partial_t+\Delta)F=0$, an integration by parts yields
		\begin{align*}
		\partial_t E[v(t)]
		&=
		-\frac14\partial_t
		\Bigl[
		\mathfrak B_4(|v+F|^2,|v+F|^2)(t)
		-\mathfrak B_4(|v|^2,|v|^2)(t)
		-\mathfrak B_4(|F|^2,|F|^2)(t)
		\Bigr]\\
		&\quad+
		\Im\int_{\mathbb R^d}
		\nabla\bigl[\overline{\mathcal H(v+F)}
		-\overline{\mathcal H(F)}\bigr]
		\cdot\nabla F\,dx.
	\end{align*}
	Consequently, by the fundamental theorem of calculus,
	\begin{align*}
		\EE(T)
		\lesssim{}&
		E[v(0)]
		+
		\Bigl\|
		\mathfrak B_4(|v+F|^2,|v+F|^2)
		-\mathfrak B_4(|v|^2,|v|^2)
		-\mathfrak B_4(|F|^2,|F|^2)
		\Bigr\|_{L_t^\infty([0,T])}\\
		&+
		\left\|
		\nabla F\cdot\nabla
		\left[
		\bigl(|x|^{-4}*|v+F|^2\bigr)(v+F)
		-\bigl(|x|^{-4}*|F|^2\bigr)F
		\right]
		\right\|_{L_{t,x}^1([0,T]\times\mathbb R^d)}.
	\end{align*}
	
	We first estimate the boundary term. Applying
	Lemma~\ref{lem:Hartree-reduction} with
	$\gamma=4$, $a=b=b'=\infty$, and
	$s=\frac{d-4}{2}$, we obtain
	\begin{align*}
		&\Bigl\|
		\mathfrak B_4(|v+F|^2,|v+F|^2)
		-\mathfrak B_4(|v|^2,|v|^2)
		-\mathfrak B_4(|F|^2,|F|^2)
		\Bigr\|_{L_t^\infty}\\
		&\quad\lesssim
		\bigl\||\nabla|^{-\frac{d-4}{2}}|v|^2
		\bigr\|_{L_t^\infty L_x^2}
		\bigl\||\nabla|^{-\frac{d-4}{2}}|F|^2
		\bigr\|_{L_t^\infty L_x^2}+
		\bigl\||\nabla|^{-\frac{d-4}{2}}|v|^2
		\bigr\|_{L_t^\infty L_x^2}
		\bigl\||\nabla|^{-\frac{d-4}{2}}
		\Re(v\overline F)\bigr\|_{L_t^\infty L_x^2}\\
		&\qquad+
		\bigl\||\nabla|^{-\frac{d-4}{2}}|F|^2
		\bigr\|_{L_t^\infty L_x^2}
		\bigl\||\nabla|^{-\frac{d-4}{2}}
		\Re(v\overline F)\bigr\|_{L_t^\infty L_x^2}+
		\bigl\||\nabla|^{-\frac{d-4}{2}}
		\Re(v\overline F)\bigr\|_{L_t^\infty L_x^2}^2\\
		&\quad\lesssim
		\EE(T)^{\frac12}
		\|F\|_{L_t^\infty L_x^{\frac{2d}{d-2}}}^2
		+
		\EE(T)^{\frac34}
		\|F\|_{L_t^\infty L_x^{\frac{2d}{d-2}}}+
		\EE(T)^{\frac14}
		\|F\|_{L_t^\infty L_x^{\frac{2d}{d-2}}}^3
		+
		\EE(T)^{\frac12}
		\|F\|_{L_t^\infty L_x^{\frac{2d}{d-2}}}^2.
	\end{align*}
	
By the Leibniz rule, it remains to estimate the following term. The first of these is the only step in the proof that requires $d\ge7$ (see
Remark~\ref{rmk:d7}). Applying \eqref{eq:hartree-grad} with $s=\tfrac{d-3}{2}$ and
$b=2$, the Sobolev embedding, and \eqref{interpolate1}, we obtain

\begin{equation}\label{5.2}
	\begin{aligned}
	\left\|
		\mathfrak B_4
		\bigl(
		\nabla|v|^2,
		(\overline v+\overline F)\nabla F
		\bigr)
		\right\|_{L_t^1}
		&\lesssim
		\bigl\||\nabla|^{-\frac{d-3}{2}}|v|^2\bigr\|_{L^2_{t,x}}\,
		\bigl\||\nabla|^{-\frac{d-7}{2}}\bigl((v+F)\nabla F\bigr)\bigr\|_{L^2_{t,x}}
		\\
		&\lesssim
		A(T)^{\frac12}
		\Bigl(
		\|v\nabla F\|_{L^2_t L_x^{\frac{2d}{2d-7}}}
		+
		\|F\nabla F\|_{L^2_t L_x^{\frac{2d}{2d-7}}}
		\Bigr)
		\\
		&\lesssim
		A(T)^{\frac12}
		\Bigl(
		\|\langle x\rangle^{-\frac3p}v\|_{L^p_tL_x^{\frac{2pd}{pd-2p+4}}}
		\|\langle x\rangle^{\frac3p}\nabla F\|_{L_t^{\frac{2p}{p-2}}L_x^{\frac{2pd}{pd-5p-4}}}
		\\
		&\hphantom{\lesssim A(T)^{\frac12}\Bigl(}
		+
		\|F\|_{L^p_tL_x^{\frac{2d}{d-2}}}
		\|\langle x\rangle^{\frac3p}\nabla F\|_{L_t^{\frac{2p}{p-2}}L_x^{\frac{2d}{d-5}}}
		\Bigr)
		\\
		&\lesssim
		\EE(T)^{\frac12-\frac1p}A(T)^{\frac12+\frac1p}
		\|F\|_{\mathcal X([0,T])}
		+
		A(T)^{\frac12}
		\|F\|_{\mathcal X([0,T])}^2 .
	\end{aligned}
\end{equation}

		For the second term, applying \eqref{eq:hartree-base} with $s=\frac{d-3}{2}$ and $b=2$, together with H\"older's inequality and \eqref{interpolate2}, we obtain
\begin{align*}
	\left\|
	\mathfrak B_4
	\bigl(
	|v|^2,
	\nabla(\overline v+\overline F)\cdot\nabla F
	\bigr)
	\right\|_{L_t^1}
	&\lesssim
	\bigl\||\nabla|^{-\frac{d-3}{2}}|v|^2\bigr\|_{L_t^2L_x^2}
	\|\x^{-\frac3p}\nabla v\|_{L_t^pL_x^2}
	\|\x^{\frac3p}\nabla F\|_{L_t^{\frac{2p}{p-2}}L_x^{\frac{2d}{d-5}}}\\
	&\quad+
	\bigl\||\nabla|^{-\frac{d-3}{2}}|v|^2\bigr\|_{L_t^{\frac p2}L_x^2}
	\|\x^{\frac3p}\nabla F\|_{L_t^{\frac{2p}{p-2}}L_x^{\frac{4d}{2d-5}}}^{2}\\
	&\lesssim
	A(T)^{\frac12}
	\EE(T)^{\frac12-\frac1p}
	A(T)^{\frac1p}
	\|F\|_{\mathcal X([0,T])}
	\\
	&\quad+
	\left\|
	|\nabla|^{-\frac{d-3}{2}}|v|^2
	\right\|_{L_t^\infty L_x^2}^{1-\frac4p}
	\left\|
	|\nabla|^{-\frac{d-3}{2}}|v|^2
	\right\|_{L_t^2L_x^2}^{\frac4p}
	\|F\|_{\mathcal X([0,T])}^{2}
	\\
	&\lesssim
	\EE(T)^{\frac12-\frac1p}A(T)^{\frac12+\frac1p}\|F\|_{\mathcal X([0,T])}
	+\EE(T)^{\frac12-\frac2p}A(T)^{\frac2p}\|F\|_{\mathcal X([0,T])}^{2}	\\
		&\quad+A(T)^{\frac2p}\|F\|_{\mathcal X([0,T])}^{2},
\end{align*}
where the second summand uses
\begin{equation}\label{eq:negative-three-energy-control-v}
	\bigl\||\nabla|^{-\frac{d-3}{2}}|v|^2\bigr\|_{L_t^\infty L_x^2}^{1-\frac4p}
	\lesssim
	\Bigl[\sup_{t\in[0,T]}\iint\frac{|v(t,x)|^2|v(t,y)|^2}{|x-y|^3}\,dx\,dy\Bigr]^{\frac12-\frac2p}
	\lesssim
	\bigl(\EE(T)+\|v\|_{L_t^\infty L_x^2}^4\bigr)^{\frac12-\frac2p}
	\lesssim
	\EE(T)^{\frac12-\frac2p}+1,
\end{equation}
obtained from the pointwise bound $|x-y|^{-3}\lesssim|x-y|^{-4}+1$ and  \eqref{mass-bd}. 

For the third,  by Lemma \ref{lem:1}, H\"older's
inequality, \eqref{interpolate1} and \eqref{interpolate2}, 
\begin{align*}
	&	\left\|
	\mathfrak B_4
	\bigl(
	\nabla|F|^2,
	\overline v\,\nabla F
	\bigr)
	+
	\mathfrak B_4
	\bigl(
	|F|^2,
	\nabla\overline v\cdot\nabla F
	\bigr)
	\right\|_{L_t^1}\\
	&\quad\lesssim
	\|F\nabla F\|_{L_t^2L_x^{\frac{d}{d-2}}}\|v\nabla F\|_{L_t^2L_x^{\frac{d}{d-2}}}
	+\||F|^2\|_{L_t^2L_x^{\frac{d}{d-2}}}\|\nabla v\,\nabla F\|_{L_t^2L_x^{\frac{d}{d-2}}}\\
	&\quad\lesssim
	\|F\|_{L_t^pL_x^{\frac{2d}{d-2}}}
	\|\x^{\frac3p}\nabla F\|_{L_t^{\frac{2p}{p-2}}L_x^{\frac{2d}{d-2}}}
	\,
	\|\x^{-\frac3p}v\|_{L_t^pL_x^{\frac{2pd}{pd-2p+4}}}
	\|\x^{\frac3p}\nabla F\|_{L_t^{\frac{2p}{p-2}}L_x^{\frac{2pd}{pd-2p-4}}}
	\\
	&\quad\quad+
	\|F\|_{L_t^4L_x^{\frac{2d}{d-2}}}^{2}
	\,
	\|\x^{-\frac3p}\nabla v\|_{L_t^pL_x^2}
	\|\x^{\frac3p}\nabla F\|_{L_t^{\frac{2p}{p-2}}L_x^{\frac{2d}{d-4}}}
	\\
	&\quad
	\lesssim
	\EE(T)^{\frac12-\frac1p}A(T)^{\frac1p}\|F\|_{\mathcal X([0,T])}^{3}.
\end{align*}

Proceeding similarly with the fourth term, by H\"older's inequality we have
\begin{equation}\label{three-energy}
	\begin{aligned}
		&	\left\|
		\mathfrak B_4
		\bigl(
		\nabla\bigl(2\Re(v\overline F)\bigr),
		(\overline v+\overline F)\nabla F
		\bigr)
		\right\|_{L_t^1}
		\\
		&\qquad\lesssim
		\left(
		\|F\nabla v\|_{L_t^2L_x^{\frac{d}{d-2}}}
		+
		\|v\nabla F\|_{L_t^2L_x^{\frac{d}{d-2}}}
		\right)
		\|(v+F)\nabla F\|_{L_t^2L_x^{\frac{d}{d-2}}}.
	\end{aligned}
\end{equation}
The first factor is bounded by
$$
\begin{aligned}
	&\|F\nabla v\|_{L_t^2L_x^{\frac{d}{d-2}}}
	+
	\|v\nabla F\|_{L_t^2L_x^{\frac{d}{d-2}}}
	\\
	&\qquad\lesssim
	\|F\|_{L_t^2L_x^{\frac{2d}{d-4}}}
	\|\nabla v\|_{L_t^\infty L_x^2}
	+
	\|\x^{-\frac3p}v\|_
	{L_t^pL_x^{\frac{2pd}{pd-2p+4}}}
	\|\x^{\frac3p}\nabla F\|_
	{L_t^{\frac{2p}{p-2}}
		L_x^{\frac{2pd}{pd-2p-4}}}
	\\
	&\qquad\lesssim
	\left(
	\EE(T)^{\frac12}
	+
	\EE(T)^{\frac12-\frac1p}A(T)^{\frac1p}
	\right)
	\|F\|_{\mathcal X([0,T])}.
\end{aligned}
$$
Similarly,
$$
\begin{aligned}
	\|(v+F)\nabla F\|_{L_t^2L_x^{\frac{d}{d-2}}}
	&\lesssim
	\|v\nabla F\|_{L_t^2L_x^{\frac{d}{d-2}}}
	+
	\|F\nabla F\|_{L_t^2L_x^{\frac{d}{d-2}}}
	\\
	&\lesssim
	\|\x^{-\frac3p}v\|_{L_t^pL_x^{\frac{2pd}{pd-2p+4}}}
	\|\x^{\frac3p}\nabla F\|_{L_t^{\frac{2p}{p-2}}L_x^{\frac{2pd}{pd-2p-4}}}
	\\
	&\quad+
	\|F\|_{L_t^pL_x^{\frac{2d}{d-2}}}
	\|\x^{\frac3p}\nabla F\|_{L_t^{\frac{2p}{p-2}}L_x^{\frac{2d}{d-2}}}
	\\
	&\lesssim
	\EE(T)^{\frac12-\frac1p}A(T)^{\frac1p}
	\|F\|_{\mathcal X([0,T])}
	+
	\|F\|_{\mathcal X([0,T])}^2.
\end{aligned}
$$
Combining the preceding estimates yields
$$
\begin{aligned}
	\mathrm{LHS}\eqref{three-energy}
	\lesssim{}&
	\EE(T)^{1-\frac1p}A(T)^{\frac1p}
	\|F\|_{\mathcal X([0,T])}^2
	+
	\EE(T)^{\frac12}
	\|F\|_{\mathcal X([0,T])}^3
	\\
	&+
	\EE(T)^{1-\frac2p}A(T)^{\frac2p}
	\|F\|_{\mathcal X([0,T])}^2
	+
	\EE(T)^{\frac12-\frac1p}A(T)^{\frac1p}
	\|F\|_{\mathcal X([0,T])}^3.
\end{aligned}
$$

Finally, for the last term,
\begin{align*}
	&	\left\|
	\mathfrak B_4
	\bigl(
	2\Re(v\overline F),
	\nabla(\overline v+\overline F)\cdot\nabla F
	\bigr)
	\right\|_{L_t^1}\\
	&\lesssim
	\|vF\|_{L_t^2L_x^{\frac{2d}{2d-3}}}
	\|\nabla v\,\nabla F\|_{L_t^2L_x^{\frac{2d}{2d-5}}}
	+
	\|vF\|_{L_t^{\frac p2}L_x^{\frac{2d}{2d-3}}}
	\|\nabla F\,\nabla F\|_{L_t^{\frac p{p-2}}L_x^{\frac{2d}{2d-5}}}
	\\
	&\lesssim
	\|v\|_{L_t^\infty L_x^2}
	\|F\|_{L_t^2L_x^{\frac{2d}{d-3}}}
	\|\x^{-\frac3p}\nabla v\|_{L_t^pL_x^2}
	\|\x^{\frac3p}\nabla F\|_{L_t^{\frac{2p}{p-2}}L_x^{\frac{2d}{d-5}}}
	\\
	&\quad+
	\|v\|_{L_t^\infty L_x^2}
	\|F\|_{L_t^{\frac p2}L_x^{\frac{2d}{d-3}}}
	\|\x^{\frac3p}\nabla F\|_{L_t^{\frac{2p}{p-2}}L_x^{\frac{4d}{2d-5}}}^{2}
	\\
	&\lesssim
	\EE(T)^{\frac12-\frac1p}
	A(T)^{\frac1p}
	\|F\|_{\mathcal X([0,T])}^{2}
	+
	\|F\|_{\mathcal X([0,T])}^{3}.
\end{align*}

Combining all the estimates above, we obtain
\begin{equation}\label{E:E1}
	\begin{aligned}
		\EE(T)
		&\lesssim
		E[v(0)]
		+\EE(T)^{\frac12}\|F\|_{L_t^\infty L_x^{\frac{2d}{d-2}}}^2
		+\EE(T)^{\frac34}\|F\|_{L_t^\infty L_x^{\frac{2d}{d-2}}}
		+\EE(T)^{\frac14}\|F\|_{L_t^\infty L_x^{\frac{2d}{d-2}}}^3
		\\
		&\quad
		+\EE(T)^{\frac12-\frac1p}A(T)^{\frac12+\frac1p}\|F\|_{\mathcal X([0,T])}
		+A(T)^{\frac12}\|F\|_{\mathcal X([0,T])}^2
		+\EE(T)^{\frac12-\frac2p}A(T)^{\frac2p}\|F\|_{\mathcal X([0,T])}^{2}
		\\
		&\quad
		+A(T)^{\frac2p}\|F\|_{\mathcal X([0,T])}^{2}
		+\EE(T)^{\frac12-\frac1p}A(T)^{\frac1p}\|F\|_{\mathcal X([0,T])}^3
		+\EE(T)^{1-\frac2p}A(T)^{\frac2p}\|F\|_{\mathcal X([0,T])}^{2}
		\\
		&\quad
	+	\EE(T)^{1-\frac1p}
		A(T)^{\frac1p}
		\|F\|_{\mathcal X([0,T])}^2
		+
		\EE(T)^{\frac12}
		\|F\|_{\mathcal X([0,T])}^3+\EE(T)^{\frac12-\frac1p}A(T)^{\frac1p}\|F\|_{\mathcal X([0,T])}^{2}
		+\|F\|_{\mathcal X([0,T])}^{3}.
	\end{aligned}
\end{equation}
By Young's inequality,  \eqref{E:E1} yields
\begin{equation}\label{E:E}
	\begin{aligned}
		\EE(T)
		&\lesssim
		E[v(0)]
		+\|F\|_{L_t^\infty L_x^{\frac{2d}{d-2}}([0,T]\times\R^d)}^{4}
		+\|F\|_{\mathcal X([0,T])}^{3}
		+A(T)^{\frac12}\|F\|_{\mathcal X([0,T])}^{2}
		+A(T)^{\frac2p}\|F\|_{\mathcal X([0,T])}^{2}
		\\
		&\quad
		+A(T)\|F\|_{\mathcal X([0,T])}^{\frac{2p}{p+2}}
		+A(T)\|F\|_{\mathcal X([0,T])}^{p}+A(T)\|F\|_{\mathcal X([0,T])}^{2p}
		+A(T)^{\frac{4}{p+4}}\|F\|_{\mathcal X([0,T])}^{\frac{4p}{p+4}}
		\\
		&\quad
		+A(T)^{\frac{2}{p+2}}\|F\|_{\mathcal X([0,T])}^{\frac{4p}{p+2}}
		+A(T)^{\frac{2}{p+2}}\|F\|_{\mathcal X([0,T])}^{\frac{6p}{p+2}}	+\|F\|_{\mathcal X([0,T])}^{6}.
	\end{aligned}
\end{equation}

	Since $F\in L_t^\infty L_x^2(\R\times\R^d)$, the mass of $u:=v+F$ is conserved, and hence
	\[
	\|v\|_{L_t^\infty L_x^2(\Imax\times\R^d)}
	\lesssim
	\|v_0\|_{L_x^2}+\|F\|_{L_t^\infty L_x^2(\R\times\R^d)}.
	\]
	Let $0<\eta\ll1$ be a small parameter. Since
	$F\in\mathcal X(\R)$, we subdivide $[0,\sup\Imax)$ into finitely many
	intervals $\{I_j\}_{j=0}^J$ such that
	\begin{equation}\label{E:partition}
		\|F\|_{\mathcal X(I_j)}\le \eta
		\qquad\text{for each }j=0,\dots,J.
	\end{equation}
	Combining \eqref{E:partition} with \eqref{E:Morawetz2}, we obtain, for $T\in I_0$,
	\begin{align}
		A(T)
		&\lesssim
		\EE(T)^{\frac12}
		+
		\EE(T)\eta^{\frac{2p}{p-2}}
		+
		\EE(T)^{\frac{p-2}{2(p-1)}}
		\left(
		\eta^{\frac{2p}{p-1}}
		+
		\eta^{\frac{3p}{p-1}}
		\right)
		\nonumber\\
		&\quad
				+
		\eta^2
		+
		\eta^4
		+
		\EE(T)
		\left(
		\eta^2+\eta^3+\eta^4
		\right).
		\label{eq:AT-bound}
	\end{align}
	Similarly, combining \eqref{E:partition} with \eqref{E:E}, we obtain
	\begin{align}
		\EE(T)
		&\lesssim
		E[v(0)]
		+
		\|F\|_{L_t^\infty L_x^{\frac{2d}{d-2}}([0,T]\times\R^d)}^{4}
		+
		\eta^3
		\nonumber\\
		&\quad
		+
		A(T)\eta^{\frac{2p}{p+2}}
		+
		A(T)\eta^{p}
		+
		A(T)^{\frac{4}{p+4}}
		\eta^{\frac{4p}{p+4}}
		\nonumber\\
		&\quad
		+
		A(T)^{\frac{2}{p+2}}
		\eta^{\frac{4p}{p+2}}
		+
		A(T)^{\frac{2}{p+2}}
		\eta^{\frac{6p}{p+2}} .
		\label{eq:ET-bound}
	\end{align}
for every \(0<T\in I_0\). In view of \eqref{3:30}, we choose $\eta>0$ sufficiently small so that a
standard continuity argument implies that $E(T)$ can grow by at most
a fixed constant on $I_0$. On each subsequent interval
$I_j=[t_j,t_{j+1}]$, the same estimates apply with $[0,T]$ replaced by
$[t_j,T]$ and $E[v(0)]$ replaced by $E[v(t_j)]$. Iterating over the
finitely many intervals $I_j$, we obtain a uniform energy bound for
$v$ on its forward maximal lifespan. The argument on the backward
maximal lifespan is identical after reversing time.
\end{proof}

\begin{remark}\label{rmk:d7}
	The threshold $d\ge 7$ in Proposition~\ref{P:energy-bds} is dictated entirely by
	the term \eqref{5.2}. The factor quadratic in $v$ is then absorbed by the interaction
	Morawetz functional
	$
	\bigl\||\nabla|^{-\frac{d-3}{2}}|v|^2\bigr\|_{L^2_{t,x}}\lesssim A(T)^{1/2}
$. Since \eqref{eq:hartree-grad} is admissible only for $0\le s\le d-5$, the choice $s=\frac{d-3}{2}$ is permissible if and only if $\frac{d-3}{2}\le d-5$, which is equivalent to $d\ge 7$.
 The remaining factor carrying $\nabla F$ then inherits the fractional smoothing $|\nabla|^{-\frac{d-7}{2}}$. By Lemma \ref{lem:1}, this operator maps $L^{\frac{2d}{2d-7}}_x$ into $L^2_x$ and is therefore available only when $\frac{d-7}{2}\ge 0$.
 
For $d\le 6$, the enforced value of $s$ lies outside the admissible range by
$\frac{d-3}{2}-(d-5)=\frac{7-d}{2}$ derivatives, corresponding to half a derivative at $d=6$
	and a full derivative at $d=5$. Consequently, the
	factor carrying $\nabla F$ would have to absorb the positive-order operator
	$|\nabla|^{\frac{7-d}{2}}$, for which no bound into $L^2_x$ is available. It remains
	open whether the dimensions $d=5,6$ can be reached by exploiting the additional
	almost sure spacetime integrability provided by the randomization of
	$F=e^{it\Delta}f^\omega$.
\end{remark}

It remains only to prove  Theorem \ref{thm:global}.
\begin{proof}[Proof of Theorem \ref{thm:global}]
	Let $d\ge 7$ and $
		1-\frac{2}{d} < s < 1$.
Let $f\in H_{rad}^s(\mathbb{R}^d)$.
	By Propositions~\ref{P:STE2} and~\ref{P:STE}, we have
	\begin{equation}
		F^\omega := e^{it\Delta}f^\omega
		\in \mathcal{X}(\mathbb{R})
		\cap L_t^\infty L_x^2(\mathbb{R}\times\mathbb{R}^d)
		\cap L_t^\infty L_x^{\frac{2d}{d-2}}(\mathbb{R}\times\mathbb{R}^d),
		\label{eq:Fomega_space}
	\end{equation}
	almost surely, where \(\mathcal{X}(\mathbb{R})\) is defined in \eqref{def:X}.	
	Fix an \(\omega\) for which \eqref{eq:Fomega_space} holds. We decompose
$
	u = F^\omega + v .
$
	Then solving \eqref{HE1} with initial data \(u(0)=f^\omega\) is equivalent to
	solving the perturbed equation \eqref{NLS-hartree-v1} for \(v\) with initial condition
	\(v(0)=0\). Proposition \ref{P:energy-bds} gives a uniform \(H_x^1\)-bound for \(v\)
	on its maximal lifespan. Consequently, Proposition \ref{P:H1-implies-scattering-hartree} implies that
	\(v\) exists globally and scatters in \(H_x^1\). This proves
	Theorem \ref{thm:global}.
\end{proof}

\section*{Acknowledgements}
The authors would like to express their sincere gratitude to Professor Zdzis{\l}aw Brze{\'z}niak for his valuable guidance and helpful discussions during a visit by the first and second authors to the University of York, where part of this work was done.

\section*{Data Availability Statement}
The data presented in this study are available in the article.


\end{document}